\documentclass[preprint,3p,12pt]{elsarticle}
\usepackage{amsmath}
\usepackage{stmaryrd}
\usepackage{bbding}
\usepackage{dcolumn}
\usepackage{graphicx}
\usepackage{amsfonts}
\usepackage{amssymb}
\usepackage{psfrag}
\usepackage{wrapfig}
\usepackage{subfigure}
\usepackage{makeidx}
\usepackage{bm}
\usepackage{epsf}
\usepackage{epsfig}
\usepackage{setspace}
\usepackage{color}

\begin{document}
\title{High-order gas-kinetic scheme in curvilinear coordinates for the Euler and Navier-Stokes solutions}

\author[BNU]{Liang Pan\corref{cor}}
\ead{panliang@bnu.edu.cn}
\address[BNU]{School of Mathematical Sciences, Beijing Normal University, Beijing, China}

\author[HKUST,HKUST2]{Kun Xu}
\ead{makxu@ust.hk}
\address[HKUST]{Department of Mathematics and Department of Mechanical and Aerospace Engineering,
The Hong Kong University of Science and Technology, Clear Water Bay, Kowloon, Hong Kong}
\address[HKUST2]{Shenzhen Research Institute, Hong Kong University of Science and Technology, Shenzhen, China}
\cortext[cor]{Corresponding author}

\begin{abstract}
The high-order gas-kinetic scheme (HGKS) has achieved success in
simulating compressible flow in Cartesian mesh. To study the flow
problem in general geometry, such as the flow over a wing-body
configuration, the development of a three-dimensional HGKS in
general curvilinear coordinates becomes necessary. In this paper, a
two-stage fourth-order gas-kinetic scheme is developed for the Euler
and Navier-Stokes solutions in the curvilinear coordinates. Based on
the coordinate transformation, the kinetic equation is transformed
first to the computational space, and the flux function in the
gas-kinetic scheme is obtained there and is transformed back to the
physical domain for the update of conservative flow variables inside
each control volume. To achieve the expected order of accuracy, the
dimension-by-dimension reconstruction based on the WENO scheme is
adopted in the computational domain, where the reconstructed
variables are the cell  averaged Jacobian and the Jacobian-weighted
conservative variables, and the conservative variables are obtained
by ratio of the above reconstructed data at Gaussian quadrature
points of each cell interface. In the two-stage fourth-order gas
kinetic scheme (GKS), similar to the generalized Riemann solver
(GRP), the initial spatial derivatives of conservative variables
have to be used in the evaluation of the time dependent flux
function in GKS, which are reconstructed as well through
orthogonalization in physical space and chain rule. A variety of
numerical examples from the order tests to the solutions with strong
discontinuities are presented to validate the accuracy and
robustness of the current scheme. The precise satisfaction of the
geometrical conservation law in non-orthogonal mesh is also
demonstrated through the numerical example.
\end{abstract}

\maketitle

\section{Introduction}
In recent decades, there have been continuous interests and efforts
on the development of high-order schemes. With the development of
computational aero-acoustics (CAA), large eddy simulations (LES), and
direct numerical simulations (DNS), the construction of high-order
numerical scheme becomes extremely demanding, and many high-order
finite volume schemes on unstructured meshes have been proposed for
the complicated geometries \cite{un-ENO,un-WENO1,un-WENO2,un-WENO3}.
However, the direct implementation in the physical space brings big
challenges. The complexity of algorithms and codes increases
dramatically because of the difficulty in choosing stencils,
especially in the multi-dimensional reconstruction. To overcome the
drawback, an efficient way is to apply the finite volume method in
the curvilinear coordinate system, where the structured meshes are
used. The technique of curvilinear or mapped coordinates is widely
used in engineering
\cite{curvilinear-1,curvilinear-2,curvilinear-3}. In principle,
given a suitable mapping function, any problem defined on a general
physical domain can be transformed into a problem in a computational
domain which is equidistant and Cartesian. Although the flexibility
may be reduced in comparison with the unstructured meshes, the good
numerical characteristics are preserved. The first one is the exact
global conservation property, which is only approximately satisfied
in the high-order finite difference method \cite{curvilinear-4}, and
the second one is the strict adherence to the integral form for
numerical simulations \cite{curvilinear-5}. Furthermore, the
standard numerical schemes on the Cartesian and equidistant grids
can be used \cite{WENO4}.

In the past decades, the gas-kinetic scheme (GKS) based on the
Bhatnagar-Gross-Krook (BGK) model \cite{BGK-1,BGK-2} has been
developed systematically for the computations from low speed flow to
supersonic one \cite{GKS-Xu1,GKS-Xu2}. Different from the numerical
methods based on the Riemann flux \cite{Riemann-appro}, GKS presents
a gas evolution process from kinetic scale to hydrodynamic scale,
where both inviscid and viscous fluxes are recovered from a
time-dependent gas distribution function at a cell interface. Based on
the unified coordinate transformation \cite{Hui}, the second-order
gas-kinetic scheme was developed under the moving-mesh framework as
well \cite{Jin1,Jin2}. The flux evaluation in the GKS is based on
the time evolution of flow variables from an initial piece-wise
discontinuous polynomials around each cell interface, where
high-order spatial and temporal evolutions of a gas distribution
function are coupled nonlinearly. With the spatial and temporal
coupled gas distribution function, the one-stage third-order GKS was
developed \cite{GKS-high-1,GKS-high-2,GKS-high-3}. In comparison
with other high-order schemes with Riemann flux
\cite{Riemann-appro}, it integrates the flux function over a time
step analytically without employing the multi-stage Runge-Kutta time
stepping techniques \cite{TVD-RK}. However, with the one-stage gas
evolution model, the formulation of GKS can become very complicated
for the further improvement, such as the one-stage fourth-order
scheme \cite{GKS-high-4}, especially for three-dimensional
computations. Based on the time-dependent flux function of the
generalized Riemann problem (GRP) \cite{GRP1,GRP2} and gas-kinetic
scheme \cite{GKS-Xu1,GKS-Xu2}, a two-stage fourth-order method was
developed for Lax-Wendroff type flow solvers \cite{Lax-Wendroff},
particularly applied for the hyperbolic conservation laws
\cite{GRP-high-1, GRP-high-2,S2O4-GKS-1}. With the temporal
discretization, a reliable framework was provided for developing GKS
into fourth-order and even higher-order accuracy with the
implementation of the traditional second-order or third-order flux
functions \cite{S2O4-GKS-2,S2O4-GKS-3,S2O4-GKS-4}. More importantly,
this scheme is as robust as the second-order scheme and works
perfectly from the subsonic to the hypersonic flows. The robustness
is due to the dynamical evolution model of the time dependent flux
function. For the construction of high-order scheme, a reliable
physical evolution model becomes important, and the delicate
flow structures captured in higher-order schemes depend on the
quality of the solvers greatly \cite{S2O4-GKS-2,S2O4-GKS-5}.

Recently, the high-order gas-kinetic scheme (HGKS) has been applied
in the direct numerical simulation of isotropic compressible
turbulence, which shows the potentials of HGKS for the simulation of
complicated flows at very high Mach numbers \cite{DNS-GKS}. To treat
practical problems with general geometry, such as the turbulent
boundary layer on non-equidistant grids and the flow over a
wing-body on non-Cartesian grids, the development of
three-dimensional HGKS in general curvilinear coordinates becomes
demanding. In this paper, based on the coordinate transformation,
the discretization procedure of finite volume method in curvilinear
coordinates is presented. To achieve the spatial accuracy, the
WENO-based dimension-by-dimension reconstruction is adopted
\cite{WENO4}, where the reconstructed variables are the cell
averaged Jacobian and the product of the conservative variables with
local Jacobian. At Gaussian quadrature points, the point value and
spatial derivatives of conservative variables weighted by local
Jacobian can be obtained. For the high-order scheme based on the
Riemann solver \cite{Riemann-appro}, the point-wise values of
conservative variables at Gaussian points are needed for the flux
calculation. However, the spatial derivatives of conservative
variables in the physical domain  is also needed in GKS, and it
plays an equally important role in the two-stage fourth-order
temporal discretization. But, it cannot be provided by the direct
spatial reconstruction.  According to the chain rule, the spatial
derivatives of conservative variables in the computational space can
be obtained first. With the procedure of orthogonalization, the
spatial derivatives in the local orthogonal coordinates of the
physical space are obtained for the flux calculation in the
second-order gas-kinetic solver. Numerical tests from
one-dimensional to three-dimensional cases, including the accuracy
test to the flows with strong discontinuities, are presented to
validate the accuracy and robustness of current scheme. The precise
satisfaction of geometrical conservation law has been demonstrated
numerically as well. The current scheme provides a solid tool for
further studies of complex compressible turbulent flows, which is
our long-term goal.

This paper is organized as follows. In Section 2, the BGK equation
and coordinate transformation are introduced. The two-stage
fourth-order gas-kinetic scheme is constructed in the curvilinear
coordinate in Section 3. Section 4 includes numerical examples to
validate the current algorithm. The last section is the
conclusion.

\section{BGK equation and coordinate transformation}
The three-dimensional BGK equation \citep{BGK-1,BGK-2} can be
written as
\begin{equation}\label{bgk}
f_t+uf_x+vf_y+wf_z=\frac{g-f}{\tau},
\end{equation}
where $\boldsymbol{u}=(u,v,w)$ is the particle velocity, $f$ is the
gas distribution function, $g$ is the three-dimensional Maxwellian
distribution and $\tau$ is the collision time. The collision term
satisfies the compatibility condition
\begin{equation}\label{compatibility}
\int \frac{g-f}{\tau}\psi \text{d}\Xi=0,
\end{equation}
where
$\displaystyle\psi=(\psi_1,...,\psi_5)^T=(1,u,v,w,\frac{1}{2}(u^2+v^2+w^2+\varsigma^2))^T$,
the internal variables
$\varsigma^2=\varsigma_1^2+...+\varsigma_K^2$,
$\text{d}\Xi=\text{d}u\text{d}v\text{d}w\text{d}\varsigma^1...\text{d}\varsigma^{K}$,
$\gamma$ is the specific heat ratio and  $K=(5-3\gamma)/(\gamma-1)$
is the degrees of freedom for three-dimensional flow. Taking
moments of the BGK equation Eq.\eqref{bgk}, the three-dimensional
conservative system can be written as
\begin{align*}
\frac{\partial Q}{\partial t}+\frac{\partial F}{\partial x}+\frac{\partial G}{\partial y}+\frac{\partial H}{\partial z}=0,
\end{align*}
where $Q=(\rho, \rho U, \rho V, \rho W,\rho E)^T$ is the
conservative variable, and $F, G, H$ are fluxes in $x, y, z$
directions given by
\begin{align*}
\begin{pmatrix}
  F \\
  G \\
  H \\
\end{pmatrix}
=\int \begin{pmatrix}
  u \\
  v \\
  w \\
\end{pmatrix}\psi f\text{d}\Xi.
\end{align*}
According to the Chapman-Enskog expansion for BGK equation, the
macroscopic governing equations can be derived
\citep{GKS-Xu1,GKS-Xu2}. In the continuum region, the BGK equation
can be rearranged and the gas distribution function can be expanded
as
\begin{align*}
f=g-\tau D_{\boldsymbol{u}}g+\tau D_{\boldsymbol{u}}(\tau
D_{\boldsymbol{u}})g-\tau D_{\boldsymbol{u}}[\tau D_{\boldsymbol{u}}(\tau
D_{\boldsymbol{u}})g]+...,
\end{align*}
where $D_{\boldsymbol{u}}=\displaystyle\frac{\partial}{\partial
t}+\boldsymbol{u}\cdot \nabla$.
With the zeroth-order truncation $f=g$, the Euler equations can ba
obtained. For the first-order truncation
\begin{align*}
f=g-\tau (ug_x+vg_y+wg_z+g_t),
\end{align*}
the Navier-Stokes equations can ba obtained. Based on the
higher-order truncations, the Burnett and super-Burnett equations
can be obtained \cite{GKS-B,GKS-SB}.

In order to construct the numerical scheme in curvilinear
coordinates, the coordinate transformation from the physical domain
$(x,y,z)$ to the computational domain $(\xi,\eta,\zeta)$ is
considered
\begin{align*}
\Big(\frac{\partial(x,y,z)}{\partial(\xi,\eta,\zeta)}\Big)=\begin{pmatrix}
x_\xi & x_\eta &x_\zeta\\
y_\xi & y_\eta &y_\zeta\\
z_\xi & z_\eta &z_\zeta\\
\end{pmatrix}.
\end{align*}
With the transformation above, the BGK equation Eq.\eqref{bgk} can
be transformed as
\begin{align}\label{bgk2}
\frac{\partial}{\partial t}(\mathcal{J} f)
+\frac{\partial}{\partial\xi}([u\widehat{\xi}_x+v\widehat{\xi}_y+w\widehat{\xi}_z]f)
&+\frac{\partial}{\partial\eta}([u\widehat{\eta}_x+v\widehat{\eta}_y+w\widehat{\eta}_z]f)\nonumber\\
&+\frac{\partial}{\partial\zeta}([u\widehat{\zeta}_x+v\widehat{\zeta}_y+w\widehat{\zeta}_z]f)
=\frac{g-f}{\tau}\mathcal{J},
\end{align}
where $\mathcal{J}$ is the Jacobian of transformation, and the
metrics above are given as follows
\begin{align*}
\begin{pmatrix}
\widehat{\xi}_x    & \widehat{\xi}_y   & \widehat{\xi}_z\\
\widehat{\eta}_x   & \widehat{\eta}_y  & \widehat{\eta}_z\\
\widehat{\zeta}_x   & \widehat{\zeta}_y  & \widehat{\zeta}_z\\
\end{pmatrix}=
\begin{pmatrix}
 y_\eta z_\zeta-z_\eta y_\zeta &  z_\eta x_\zeta-x_\eta z_\zeta &  x_\eta y_\zeta-y_\eta x_\zeta \\
 y_\zeta z_\xi-z_\zeta y_\xi &  z_\zeta x_\xi-x_\zeta z_\xi &  x_\zeta y_\xi-y_\zeta x_\xi \\
 y_\xi z_\eta-z_\xi y_\eta &  z_\xi x_\eta-x_\xi z_\eta &  x_\xi y_\eta-y_\xi x_\eta \\
\end{pmatrix}.
\end{align*}
Taking moments of Eq.\eqref{bgk2}, the macroscopic equations can be
written as
\begin{align}\label{computation}
\frac{\partial \widehat{Q}}{\partial t}+\frac{\partial
\widehat{F}}{\partial \xi}+\frac{\partial \widehat{G}}{\partial
\eta}+\frac{\partial \widehat{H}}{\partial \zeta}=0,
\end{align}
where $\widehat{Q}=\mathcal{J}Q$ and $\widehat{F}, \widehat{G},
\widehat{H}$ are fluxes in $\xi,\eta, \zeta$ directions given by
\begin{align*}
\widehat{F}=\int[u\widehat{\xi}_x+v\widehat{\xi}_y+w\widehat{\xi}_z]\psi f\text{d}\Xi,\\
\widehat{G}=\int[u\widehat{\eta}_x+v\widehat{\eta}_y+w\widehat{\eta}_z]\psi f\text{d}\Xi, \\
\widehat{H}=\int[u\widehat{\zeta}_x+v\widehat{\zeta}_y+w\widehat{\zeta}_z]\psi f\text{d}\Xi.
\end{align*}

\begin{figure}[!h]
\centering
\includegraphics[width=0.8\textwidth]{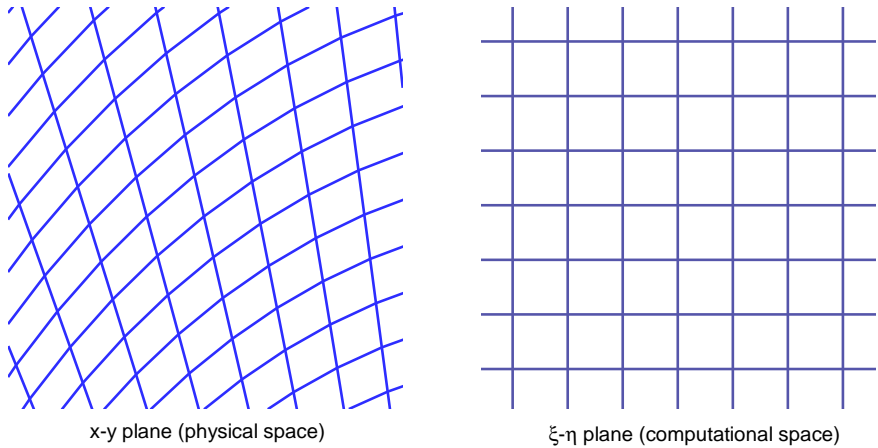}
\caption{\label{transformation} Schematic of physical domain $(x,y)$
and computational domain $(\xi,\eta)$ in two-dimensional case.}
\end{figure}

Integrating Eq.\eqref{computation} over the control volume
$V_{ijk}$, the semi-discretized finite volume scheme can be written
as
\begin{align}\label{finite}
\frac{\text{d}\widehat{Q}_{ijk}}{\text{d}t}=\mathcal{L}(\widehat{Q}_{ijk})=-\frac{1}{|V_{ijk}|}\Big[
&\int_{\eta_j-\Delta \eta/2}^{\eta_j+\Delta \eta/2}\int_{\zeta_k-\Delta \zeta/2}^{\zeta_k+\Delta \zeta/2}(\widehat{F}_{i+1/2,j,k}-\widehat{F}_{i-1/2,j,k})\text{d}\eta\text{d}\zeta\nonumber\\
+&\int_{\xi_i-\Delta \xi/2}^{\xi_i+\Delta \xi/2}\int_{\zeta_k-\Delta \zeta/2}^{\zeta_k+\Delta \zeta/2}(\widehat{G}_{i,j+1/2,k}-\widehat{G}_{i,j-1/2,k})\text{d}\xi\text{d}\zeta\nonumber\\
+&\int_{\xi_i-\Delta \xi/2}^{\xi_i+\Delta \xi/2}\int_{\eta_j-\Delta \eta/2}^{\eta_j+\Delta \eta/2}(\widehat{H}_{i,j,k+1/2}-\widehat{H}_{i,j,k-1/2})\text{d}\xi\text{d}\eta\Big],
\end{align}
where $|V_{ijk}|=\Delta \xi\Delta \eta\Delta \zeta$. The numerical
fluxes in $\xi$-direction is given as example. To achieve the
spatial accuracy, the Gaussian quadrature is used over the cell
interface and we have
\begin{align}\label{flux-x}
\int_{\eta_j-\Delta \eta/2}^{\eta_j+\Delta
\eta/2}\int_{\zeta_k-\Delta \zeta/2}^{\zeta_k+\Delta
\zeta/2}\widehat{F}_{i+1/2,j,k}\text{d}\eta\text{d}\zeta=\Delta\eta\Delta\zeta\sum_{m,n=1}^2\omega_{mn}
\widehat{F}(\boldsymbol{\xi}_{i+1/2,j_m,k_n},t),
\end{align}
where $\boldsymbol{\xi}_{i+1/2,j_m,k_n}$ is the Gauss quadrature
point of cell interface $[\eta_j-\Delta \eta/2,\eta_j+\Delta
\eta/2]\times[\zeta_k-\Delta \zeta/2,\zeta_k+\Delta \zeta/2]$ with
$\xi=\xi_{i+1/2}$, and $\omega_{mn}$ are quadrature weights.
According to the definition of $\widehat{F}$, the numerical flux in
Eq.\eqref{flux-x} for each quadrature point can rewritten as
\begin{align*}
\widehat{F}(\boldsymbol{\xi}_{i+1/2,j_m,k_n},t)=S_{\xi}\int
\widetilde{u}\psi
f(\boldsymbol{x}_{i+1/2,j_m,k_n},t,\widetilde{\boldsymbol{u}},\varsigma)\text{d}\widetilde{\Xi},
\end{align*}
where
$S_{\xi}=\sqrt{\widehat{\xi}_x^2+\widehat{\xi}_y^2+\widehat{\xi}_z^2}$
and the local particle velocity can be given by
\begin{align*}
(\widetilde{u},\widetilde{v},\widetilde{w})=(u,v,w)\cdot(\boldsymbol{n}_x,\boldsymbol{n}_y,
\boldsymbol{n}_z),
\end{align*}
where $\boldsymbol{n}_x$ is the normal direction, and
$\boldsymbol{n}_y, \boldsymbol{n}_z$ are two orthogonal tangential
directions at each Gaussian quadrature point, which can be
determined sequentially
\begin{align*}
\boldsymbol{n}_x&=(\widehat{\xi}_x,\widehat{\xi}_y,\widehat{\xi}_z)/\sqrt{\widehat{\xi}_x^2+\widehat{\xi}_y^2+\widehat{\xi}_z^2},\\
\boldsymbol{n}_z&=(x_\zeta, y_\zeta, z_\zeta)/\sqrt{x_\zeta^2+y_\zeta^2+z_\zeta^2},\\
\boldsymbol{n}_y&=\boldsymbol{n}_z\times\boldsymbol{n}_x.
\end{align*}
Denote $(a_{ij})$ is the inverse of
$(\boldsymbol{n}_x,\boldsymbol{n}_y, \boldsymbol{n}_z)$, and each
component of $\widehat{F}(\boldsymbol{\xi}_{i+1/2,j_m,k_n},t)$ can
be given by the combination of fluxes in the local orthogonal
coordinate
\begin{align*}
\left\{\begin{aligned}
F_{\rho}&=S_{\xi}F_{\widetilde{\rho}},\\
F_{\rho u}&=S_{\xi}(a_{11}F_{\widetilde{\rho u}}+a_{12}F_{\widetilde{\rho v}}+a_{13}F_{\widetilde{\rho w}}),\\
F_{\rho v}&=S_{\xi}(a_{21}F_{\widetilde{\rho v}}+a_{22}F_{\widetilde{\rho u}}+a_{23}F_{\widetilde{\rho w}}),\\
F_{\rho v}&=S_{\xi}(a_{31}F_{\widetilde{\rho v}}+a_{32}F_{\widetilde{\rho u}}+a_{33}F_{\widetilde{\rho w}}),\\
F_{E}&=S_{\xi}F_{\widetilde{E}},
\end{aligned} \right.
\end{align*}
where the fluxes in the local coordinate can be obtained as follows
\begin{align*}
(F_{\widetilde{\rho}},
    F_{\widetilde{\rho u}},
    F_{\widetilde{\rho v}},
    F_{\widetilde{\rho w}},
    F_{\widetilde{E}})^T=\int\widetilde{u}\widetilde{\psi} f(\boldsymbol{x}_{i+1/2,j_m,k_n},t,\widetilde{\boldsymbol{u}},\varsigma)\text{d}\widetilde{\Xi},
\end{align*}
and
$\widetilde{\psi}=(1,\widetilde{u},\widetilde{v},\widetilde{w},(\widetilde{u}^2+\widetilde{v}^2+\widetilde{w}^2+\varsigma^2)/2)^T$.
The procedure above shows that the spatial reconstruction, including
the conservative variables and their spatial derivatives, needs to
be conducted in the orthogonal local coordinate $(\boldsymbol{n}_x,
\boldsymbol{n}_y, \boldsymbol{n}_z)$ in the physical domain. With
the integral solution of BGK equation, the gas distribution function
can be constructed as follows
\begin{equation*}
f(\boldsymbol{x}_{i+1/2,j_m,k_n},t,\boldsymbol{u},\varsigma)=\frac{1}{\tau}\int_0^t
g(\boldsymbol{x}',t',\boldsymbol{u},\varsigma)e^{-(t-t')/\tau}\text{d}t'+e^{-t/\tau}f_0(-\boldsymbol{u}t,\varsigma),
\end{equation*}
where
$\widetilde{\boldsymbol{u}}=(\widetilde{u},\widetilde{v},\widetilde{w})$
is denoted as $\boldsymbol{u}=(u,v,w)$ for simplicity in this
section,
$\boldsymbol{x}_{i+1/2,j_m,k_n}=(x_{i+1/2},y_{j_m},z_{k_n})$ is the
location of Gaussian quadrature point, $x_{i+1/2}=x'+u(t-t'),
y_{j_m}=y'+v(t-t'), z_{k_n}=z'+w(t-t')$ are the trajectory of
particles, $f_0$ is the initial gas distribution function, and $g$
is the corresponding equilibrium state. With the reconstruction of
macroscopic variables, the gas distribution function at the cell
interface can be expressed as
\begin{align}\label{flux}
f(\boldsymbol{x}_{i+1/2,j_m,k_n},t,\boldsymbol{u},\varsigma)=&(1-e^{-t/\tau})g_0+((t+\tau)e^{-t/\tau}-\tau)(\overline{a}_1u+\overline{a}_2v+\overline{a}_3w)g_0\nonumber\\
+&(t-\tau+\tau e^{-t/\tau}){\bar{A}} g_0\nonumber\\
+&e^{-t/\tau}g_r[1-(\tau+t)(a_{1}^{r}u+a_{2}^{r}v+a_{3}^{r}w)-\tau A^r)]H(u)\nonumber\\
+&e^{-t/\tau}g_l[1-(\tau+t)(a_{1}^{l}u+a_{2}^{l}v+a_{3}^{l}w)-\tau
A^l)](1-H(u)),
\end{align}
where the equilibrium state $g_{0}$ and corresponding conservative
variables $Q_{0}$ at the quadrature point can be determined by the
compatibility condition Eq.\eqref{compatibility}
\begin{align*}
\int\psi g_{0}\text{d}\Xi=Q_0=\int_{u>0}\psi
g_{l}\text{d}\Xi+\int_{u<0}\psi g_{r}\text{d}\Xi,
\end{align*}
and the coefficients in Eq.\eqref{flux} can be determined by the
reconstructed directional derivatives and compatibility condition
\begin{align*}
\displaystyle \langle a_{1}^{k}\rangle=\frac{\partial
Q_{k}}{\partial \boldsymbol{n}_x}, \langle
a_{2}^{k}\rangle=\frac{\partial Q_{k}}{\partial \boldsymbol{n}_y},
\langle a_{3}^{k}\rangle&=\frac{\partial Q_{k}}{\partial
\boldsymbol{n}_z}, \langle
a_{1}^{k}u+a_{2}^{k}v+a_{3}^{k}w+A^{k}\rangle=0,\\ \displaystyle
\langle\overline{a}_1\rangle=\frac{\partial Q_{0}}{\partial
{\boldsymbol{n}_x}}, \langle\overline{a}_2\rangle=\frac{\partial
Q_{0}}{\partial {\boldsymbol{n}_y}},
\langle\overline{a}_3\rangle&=\frac{\partial Q_{0}}{\partial
{\boldsymbol{n}_z}},
\langle\overline{a}_1u+\overline{a}_2v+\overline{a}_3w+\overline{A}\rangle=0,
\end{align*}
where $k=l,r$ and $\langle...\rangle$ are the moments of the
equilibrium $g$ and defined by
\begin{align*}
\langle...\rangle=\int g (...)\psi \text{d}\Xi.
\end{align*}
More details of the gas-kinetic scheme can be found in
\cite{GKS-Xu1}.

\section{High-order scheme in curvilinear coordinate}
\subsection{Spatial reconstruction}
To achieve the high-order spatial accuracy, the fifth-order WENO
reconstruction \cite{WENO1, WENO2, WENO3} is adopted. In the
curvilinear coordinates, the reconstruction is conducted for the
cell averaged variables $\mathcal{J} Q$ and $\mathcal{J}$. However,
the reconstructed variables $Q$ and spatial derivatives
$\displaystyle\frac{\partial Q}{\partial
\boldsymbol{n}_x},\frac{\partial Q}{\partial
\boldsymbol{n}_y},\frac{\partial Q}{\partial \boldsymbol{n}_z}$ are
needed for the gas-kinetic solver, and the special treatment is
needed for reconstruction. The procedure is given as follows
\begin{enumerate}
\item For each Gaussian quadrature point
$\boldsymbol{\xi}_{i+1/2,j_m,k_n}=(\xi_{i+1/2},\eta_{j_m},\zeta_{k_n})$,
the local coordinate $(\boldsymbol{n}_x,\boldsymbol{n}_y,
\boldsymbol{n}_z)$ is determined. For a general coordinate
transformation, the local coordinate is different for each
quadrature point. More computational cost will be introduced for the
reconstruction because the variables need to be projected into
different local coordinate.
\item According to one-dimensional WENO-Z reconstruction \cite{WENO3}, the
cell averaged reconstructed values and cell averaged spatial
derivatives at $\xi=\xi_{i+1/2}$ can be constructed
\begin{align*}
(\mathcal{J}& Q_{l})_{j-\ell_1,k-\ell_2},
(\mathcal{J} Q_{r})_{j-\ell_1,k-\ell_2},(\mathcal{J} Q_{0})_{j-\ell_1,k-\ell_2},\\
(\partial_\xi(\mathcal{J}Q)_{l}&)_{j-\ell_1,k-\ell_2},
(\partial_\xi(\mathcal{J}Q)_{r})_{j-\ell_1,k-\ell_2},(\partial_\xi(\mathcal{J}Q)_{0})_{j-\ell_1,k-\ell_2},
\end{align*}
where $\ell_1,\ell_2=-2,...,2$. With the WENO reconstruction in the
horizontal direction over the interval $[\zeta_{k-\ell_2}-\Delta
\zeta/2, \zeta_{\zeta-\ell_2}+\Delta \zeta/2]$, the averaged value
and the averaged spatial derivatives with $\eta=\eta_{j_m}$ can be
given
\begin{align*}
(\mathcal{J}Q_{l})_{j_m,k-\ell_2},(&\mathcal{J}Q_{r})_{j_m,k-\ell_2},(\mathcal{J}Q_{0})_{j_m,k-\ell_2},\\
(\partial_\xi(\mathcal{J}Q)_{l})_{j_m,k-\ell_2},(&\partial_\xi(\mathcal{J}Q)_{r})_{j_m,k-\ell_2},(\partial_\xi(\mathcal{J}Q)_{0})_{j_m,k-\ell_2},\\
(\partial_\eta(\mathcal{J}Q)_{l})_{j_m,k-\ell_2},(&\partial_\eta(\mathcal{J}Q)_{r})_{j_m,k-\ell_2},(\partial_\eta(\mathcal{J}Q)_{0})_{j_m,k-\ell_2}.
\end{align*}
With the WENO reconstruction in the vertical direction, the point
value and spatial derivatives at Gaussian quadrature points
$\boldsymbol{\xi}_{i+1/2,m,n}=(\xi_{i+1/2},\eta_{j_m},\zeta_{k_n})$
can be given
\begin{align*}
(\mathcal{J}Q_{l})_{j_m,k_n}, &(\mathcal{J}Q_{r})_{j_m,k_n},(\mathcal{J}Q_{0})_{j_m,k_n},\\
(\partial_\xi(\mathcal{J}Q)_{l})_{j_m,k_n},&(\partial_\xi(\mathcal{J}Q)_{r})_{j_m,k_n},(\partial_\xi(\mathcal{J}Q)_{0})_{j_m,k_n},\\
(\partial_\eta(\mathcal{J}Q)_{l})_{j_m,k_n},&(\partial_\eta(\mathcal{J}Q)_{r})_{j_m,k_n},(\partial_\eta(\mathcal{J}Q)_{0})_{j_m,k_n},\\
(\partial_\zeta(\mathcal{J}Q)_{l})_{j_m,k_n},&(\partial_\zeta(\mathcal{J}Q)_{r})_{j_m,k_n},(\partial_\zeta(\mathcal{J}Q)_{0})_{j_m,k_n}.
\end{align*}
More details of three-dimensional high-order gas-kinetic scheme can
be found in \cite{GKS-high-2}.
\item With the identical procedure, the reconstructed Jacobian can be
obtained at the Gaussian quadrature point
$\boldsymbol{\xi}_{i+1/2,j_m,k_n}=(\xi_{i+1/2},\eta_{j_m},\zeta_{k_n})$
as well
\begin{align*}
(\mathcal{J}_{l})_{j_m,k_n}, &(\mathcal{J}_{r})_{j_m,k_n},(\mathcal{J}_{0})_{j_m,k_n},\\
(\partial_\xi\mathcal{J}_{l})_{j_m,k_n},&(\partial_\xi\mathcal{J}_{r})_{j_m,k_n},(\partial_\xi\mathcal{J}_{0})_{j_m,k_n},\\
(\partial_\eta\mathcal{J}_{l})_{j_m,k_n},&(\partial_\eta\mathcal{J}_{r})_{j_m,k_n},(\partial_\eta\mathcal{J}_{0})_{j_m,k_n},\\
(\partial_\zeta\mathcal{J}_{l})_{j_m,k_n},&(\partial_\zeta\mathcal{J}_{r})_{j_m,k_n},(\partial_\zeta\mathcal{J}_{0})_{j_m,k_n}.
\end{align*}
\item For simplicity, the subscripts corresponding to the Gaussian
quadrature points and the variables at the left, right and across
the cell interface are omitted. With the reconstruction of
$(\mathcal{J} Q)$ and $\mathcal{J}$, the point value $Q$ can be
calculated by
\begin{align*}
Q=\frac{(\mathcal{J} Q)}{\mathcal{J}}.
\end{align*}
The spatial derivatives $Q_\xi, Q_\eta, Q_\zeta$ in the
computational domain can be obtained by the above reconstructed
$(\mathcal{J} Q), \mathcal{J}$, and chain rule as well
\begin{align*}
Q_\xi=\frac{(\mathcal{J} Q)_{\xi}-Q\mathcal{J}_\xi}{\mathcal{J}},\\
Q_\eta=\frac{(\mathcal{J} Q)_{\eta}-Q\mathcal{J}_\eta}{\mathcal{J}},\\
Q_\zeta=\frac{(\mathcal{J}
Q)_{\eta}-Q\mathcal{J}_\zeta}{\mathcal{J}}.
\end{align*}
However, what we need is the directional derivatives
$\displaystyle\frac{\partial Q}{\partial
\boldsymbol{n}_x},\frac{\partial Q}{\partial
\boldsymbol{n}_y},\frac{\partial Q}{\partial \boldsymbol{n}_z}$.
According to the chain rule, the spatial derivatives can be
rewritten as
\begin{align*}
Q_\xi&=Q_x x_\xi+Q_y y_\xi+Q_z z_\xi,\\
Q_\eta&=Q_x x_\eta+Q_y y_\eta+Q_z z_\eta,\\
Q_\zeta&=Q_x x_\zeta+Q_y y_\zeta+Q_z z_\zeta.
\end{align*}
The normalized spatial derivatives can be considered as the
directional derivatives along the following direction
\begin{align*}
Q_{\xi'}=Q_\xi/|\boldsymbol{x}_\xi|,~~&\boldsymbol{\tau}_1=(x_\xi, y_\xi,z_\xi)/|\boldsymbol{x}_\xi|\\
Q_{\eta'}=Q_\eta/|\boldsymbol{x}_\eta|,~~&\boldsymbol{\tau}_2=(x_\eta, y_\eta,z_\eta)/|\boldsymbol{x}_\eta|\\
Q_{\zeta'}=Q_\zeta/|\boldsymbol{x}_\zeta|,~~&\boldsymbol{\tau}_3=(x_\zeta,
y_\zeta, z_\zeta)/|\boldsymbol{x}_\zeta|,
\end{align*}
where $\boldsymbol{\tau}_1, \boldsymbol{\tau}_2,
\boldsymbol{\tau}_3$ can be obtained from the coordinate
transformation. For the Cartesian mesh, they coincide with
$\boldsymbol{n}_x,\boldsymbol{n}_y, \boldsymbol{n}_z$. However, for
the general meshes, they are not orthogonal. The procedure of
orthogonalization is used to generate the spatial derivatives in the
local orthogonal coordinate for the calculation of numerical fluxes
\begin{align*}
\frac{\partial Q}{\partial
\boldsymbol{n}_z}&=Q_{\zeta'},\\
\frac{\partial Q}{\partial
\boldsymbol{n}_y}&=\frac{1}{(\boldsymbol{\tau}_2,\boldsymbol{n}_y)}Q_{\eta'}-\frac{(\boldsymbol{\tau}_2,\boldsymbol{n}_z)}{(\boldsymbol{\tau}_2,\boldsymbol{n}_y)}\frac{\partial
Q}{\partial
\boldsymbol{n}_z},\\
\frac{\partial Q}{\partial
\boldsymbol{n}_x}&=\frac{1}{(\boldsymbol{\tau}_1,\boldsymbol{n}_x)}Q_{\xi'}-\frac{(\boldsymbol{\tau}_1,\boldsymbol{n}_y)}{(\boldsymbol{\tau}_1,\boldsymbol{n}_x)}\frac{\partial
Q}{\partial
\boldsymbol{n}_y}-\frac{(\boldsymbol{\tau}_1,\boldsymbol{n}_z)}{(\boldsymbol{\tau}_1,\boldsymbol{n}_x)}\frac{\partial
Q}{\partial \boldsymbol{n}_z}.
\end{align*}
Thus, the spatial derivatives in the local orthogonal coordinate are
fully determined. The fourth step is analytical and no error is
introduced. So long as the spatial accuracy is achieved in the
second and third step, the order of accuracy can be maintained by the
procedures above.
\end{enumerate}

\subsection{Temporal discretization}
A two-stage fourth-order time-accurate discretization was developed
for  Lax-Wendroff flow solvers with the generalized Riemann problem
(GRP) solver \cite{GRP-high-1} and the gas-kinetic scheme (GKS)
\cite{S2O4-GKS-1}. Consider the following time-dependent equation
\begin{align*}
\frac{\partial \boldsymbol{q}}{\partial t}=\mathcal {L}(\boldsymbol{q}),
\end{align*}
with the initial condition at $t_n$, i.e.,
\begin{align*}
\boldsymbol{q}(t=t_n)=\boldsymbol{q}^n,
\end{align*}
where $\mathcal {L}$ is an operator for spatial derivative of flux.
The state $\boldsymbol{q}^{n+1}$ at $t_{n+1}=t_n+\Delta t$  can be
updated with the following formula
\begin{align*}
&\boldsymbol{q}^*=\boldsymbol{q}^n+\frac{1}{2}\Delta t\mathcal
{L}(\boldsymbol{q}^n)+\frac{1}{8}\Delta t^2\frac{\partial}{\partial
t}\mathcal{L}(\boldsymbol{q}^n), \\
\boldsymbol{q}^{n+1}=&\boldsymbol{q}^n+\Delta t\mathcal
{L}(\boldsymbol{q}^n)+\frac{1}{6}\Delta t^2\big(\frac{\partial}{\partial
t}\mathcal{L}(\boldsymbol{q}^n)+2\frac{\partial}{\partial
t}\mathcal{L}(\boldsymbol{q}^*)\big).
\end{align*}
It can be proved that for hyperbolic equations the above temporal
discretization provides a fourth-order time accurate solution for
$\boldsymbol{q}^{n+1}$.

In order to develop the high-order scheme in the curvilinear
coordinate, the semi-discretized finite volume scheme
\begin{align*}
\frac{\text{d}\widehat{Q}_{ijk}}{\text{d}t}=\mathcal{L}(\widehat{Q}_{ijk}),
\end{align*}
can be discretized according to the two-stage temporal method.  To
implement two-stage method, the following notation is introduced
\begin{align*}
\widehat{\mathbb{F}}(\boldsymbol{\xi}_{i+1/2,j,k},\delta)
=\sum_{m,n=1}^2 S_{mn}\int_{t_n}^{t_n+\delta}\int \widetilde{u}\psi
f(\boldsymbol{x}_{i+1/2,j_m,k_n},t,\widetilde{\boldsymbol{u}},\varsigma)\text{d}\widetilde{\Xi}\text{d}t.
\end{align*}
and it can be expanded as the following linear form in the time
interval $[t_n, t_n+\Delta t]$
\begin{align*}
\widehat{\mathbb{F}}(\boldsymbol{\xi}_{i+1/2,j,k},t)=\widehat{F}_{i+1/2,j,k}^n+ \partial_t \widehat{F}_{i+1/2,j,k}^n(t-t_n).
\end{align*} Integrate over  $[t_n, t_n+\Delta t/2]$ and  $[t_n, t_n+\Delta
t]$, we have the following two equations
\begin{align*}
\widehat{F}_{i+1/2,j,k}^n\Delta t&+\frac{1}{2}\partial_t \widehat{F}_{i+1/2,j,k}^n\Delta t^2 =\widehat{\mathbb{F}}(\boldsymbol{\xi}_{i+1/2,j,k},\Delta t) , \\
\frac{1}{2}\widehat{F}_{i+1/2,j,k}^n\Delta t&+\frac{1}{8}\partial_t \widehat{F}_{i+1/2,j,k}^n\Delta t^2 =\widehat{\mathbb{F}}(\boldsymbol{\xi}_{i+1/2,j,k}, \Delta t/2).
\end{align*}
By solving the linear system, The coefficients
$\widehat{F}_{i+1/2,j,k}^n$ and $\partial_t
\widehat{F}_{i+1/2,j,k}^n$ can be determined. Similarly,
$\widehat{F}_{i+1/2,j,k}^*$ and $\partial_t
\widehat{F}_{i+1/2,j,k}^*$ for the intermediate state can be
constructed. More details of the two-stage fourth-order scheme can
be found in \cite{GRP-high-1, S2O4-GKS-1}

\subsection{One-dimensional scheme in non-equidistant grids}
As a particular case, the method for one-dimensional flow
degenerates to the scheme in non-equidistant grids. For
one-dimensional flows, the finite volume scheme Eq.\eqref{finite}
can be simplified as
\begin{align*}
\frac{\text{d}(\mathcal{J}Q)_{i}}{\text{d}
t}=-\frac{1}{\Delta\xi}(F_{i+1/2}-F_{i-1/2}),
\end{align*}
where $F_{i+1/2}$ is the numerical flux in the physical domain. To
implement the high-order gas-kinetic scheme, the point value
$Q_{i+1/2}$ and spatial derivative $(Q_x)_{i+1/2}$ are needed. With
the reconstructed for cell averaged variable $(\mathcal{J} Q)$ and
cell averaged Jacobian $\mathcal{J}$, the point value is given by
\begin{align*}
Q_{i+1/2}&=\frac{(\mathcal{J} Q)_{i+1/2}}{\mathcal{J}_{i+1/2}}.
\end{align*}
With the coordinate transformation, the relation of spatial
derivative can be expressed as
\begin{align*}
Q_\xi=Q_xx_\xi.
\end{align*}
The spatial derivative can be calculated by
\begin{align*}
(Q_x)_{i+1/2}=\frac{(Q_\xi)_{i+1/2}}{(x_\xi)_{i+1/2}},
\end{align*}
where $(Q_\xi)_{i+1/2}$ is given by the chain rule
\begin{align*}
(Q_\xi)_{i+1/2}=\frac{((\mathcal{J}
Q)_{\xi})_{i+1/2}-Q_{i+1/2}(\mathcal{J}_\xi)_{i+1/2}}{\mathcal{J}_{i+1/2}}.
\end{align*}
With the above procedure, the one-dimensional gas-kinetic scheme is obtained.

\section{Numerical tests}
In this section, numerical tests for both inviscid and viscous
flows will be presented to validate our numerical scheme. For the
inviscid flow, the collision time $\tau$ takes
\begin{align*}
\tau=\epsilon \Delta t+C\displaystyle|\frac{p_l-p_r}{p_l+p_r}|\Delta
t,
\end{align*}
where $\varepsilon=0.01$ and $C=1$. For the viscous flow, we have
\begin{align*}
\tau=\frac{\nu}{p}+C \displaystyle|\frac{p_l-p_r}{p_l+p_r}|\Delta t,
\end{align*}
where $p_l$ and $p_r$ denote the pressure on the left and right
sides of the cell interface, $\nu$ is the dynamic viscous
coefficient, and $p$ is the pressure at the cell interface. The
ratio of specific heats takes $\gamma=1.4$. The reason for including
artificial dissipation through the additional term in the particle
collision time is to enlarge the kinetic scale physics in the
discontinuous region for the construction of a numerical shock
structure through the particle free transport and inadequate
particle collision in order to keep the non-equilibrium property.

\begin{table}[!h]
\begin{center}
\def\temptablewidth{0.75\textwidth}
{\rule{\temptablewidth}{1.0pt}}
\begin{tabular*}{\temptablewidth}{@{\extracolsep{\fill}}c|cc|cc}
mesh & $L^1$ error & order ~ & $L^2$ error & order        \\
\hline
10  &  2.5450E-03 &  ~~      &   2.0040E-03  &   ~~      \\
20  &  8.0378E-05 &  4.9847  &   6.3372E-05  &   4.9829  \\
40  &  2.5856E-06 &  4.9582  &   2.0277E-06  &   4.9659  \\
80  &  8.1489E-08 &  4.9877  &   6.3762E-08  &   4.9910  \\
160 &  2.5499E-09 &  4.9980  &   1.9959E-09  &   4.9975  \\
320 &  7.9780E-11 &  4.9983  &   6.2447E-11  &   4.9982  \\
640 &  2.5150E-12 &  4.9873  &   1.9801E-12  &   4.9789  \\
\end{tabular*}
{\rule{\temptablewidth}{1.0pt}}
\end{center}
\vspace{-5mm}\caption{\label{tab-1d-1} Accuracy test: 1D advection
of density perturbation with nonuniform meshes.}
\begin{center}
\def\temptablewidth{0.75\textwidth}
{\rule{\temptablewidth}{1.0pt}}
\begin{tabular*}{\temptablewidth}{@{\extracolsep{\fill}}c|cc|cc}
mesh & $L^1$ error & order ~ & $L^2$ error & order      \\
\hline
10  &   1.7066E-03 &   ~~     &   1.3691E-03  &  ~~     \\
20  &   5.7014E-05 &  4.9036  &   4.4846E-05  & 4.9321  \\
40  &   1.8059E-06 &  4.9804  &   1.4147E-06  & 4.9863  \\
80  &   5.6518E-08 &  4.9979  &   4.4293E-08  & 4.9972  \\
160 &   1.7678E-09 &  4.9986  &   1.3852E-09  & 4.9989  \\
320 &   5.5340E-11 &  4.9975  &   4.3366E-11  & 4.9973  \\
640 &   1.7396E-12 &  4.9914  &   1.3633E-12  & 4.9913  \\
\end{tabular*}
{\rule{\temptablewidth}{1.0pt}}
\end{center}
\vspace{-5mm}\caption{\label{tab-1d-2} Accuracy test: 1D advection
of density perturbation with uniform meshes.}
\end{table}

\subsection{Accuracy tests}
The advection of density perturbation for the one-dimensional to
three-dimensional flows are presented to test the order of accuracy.
For the one-dimensional case, the physical domain is $[0,2]$ and the
initial conditions are set as follows
\begin{align*}
\rho_0(x)=1+0.2\sin(\pi x),~p_0(x)=1,~U_0(x)=1.
\end{align*}
The periodic boundary conditions are imposed at both ends of the
physical domain and the exact solutions are
\begin{align*}
\rho(x,t)=1+0.2\sin(\pi(x-t)),~p(x,t)=1,~U(x,t)=1.
\end{align*}
The computational domain is $[0,2]$ as well, and a nonuniform mesh
is provided by the following coordinate transformation
\begin{align*}
x=\xi+0.05\sin(\pi\xi),
\end{align*}
where $N$ uniform cells are used in computational domain. In order
to get the cell integrated flow variables, three-point Gaussian
quadrature is used inside each cell to evaluate the values without
losing accuracy. As reference, the mesh with $N$ uniform cells in
physical domain is tested as well. The $L^1$ and $L^2$ errors and
orders of accuracy at $t=2$ are presented in Tab.\ref{tab-1d-1} and
Tab.\ref{tab-1d-2} for both nonuniform and uniform meshes,
respectively. The expected order of accuracy are achieved with the
mesh refinement.

\begin{figure}[!htb]
\centering
\includegraphics[width=0.48\textwidth]{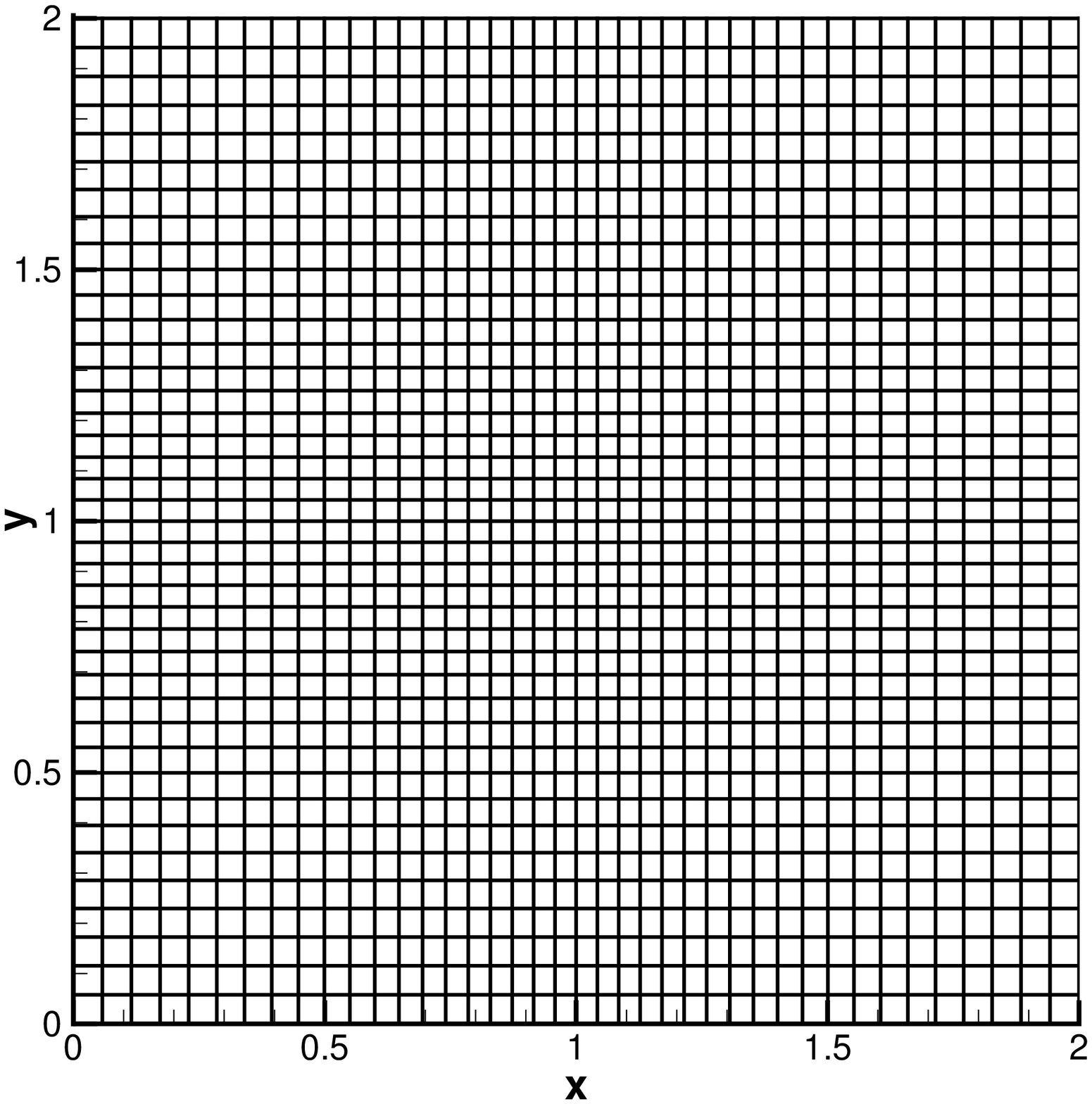}{a}
\includegraphics[width=0.48\textwidth]{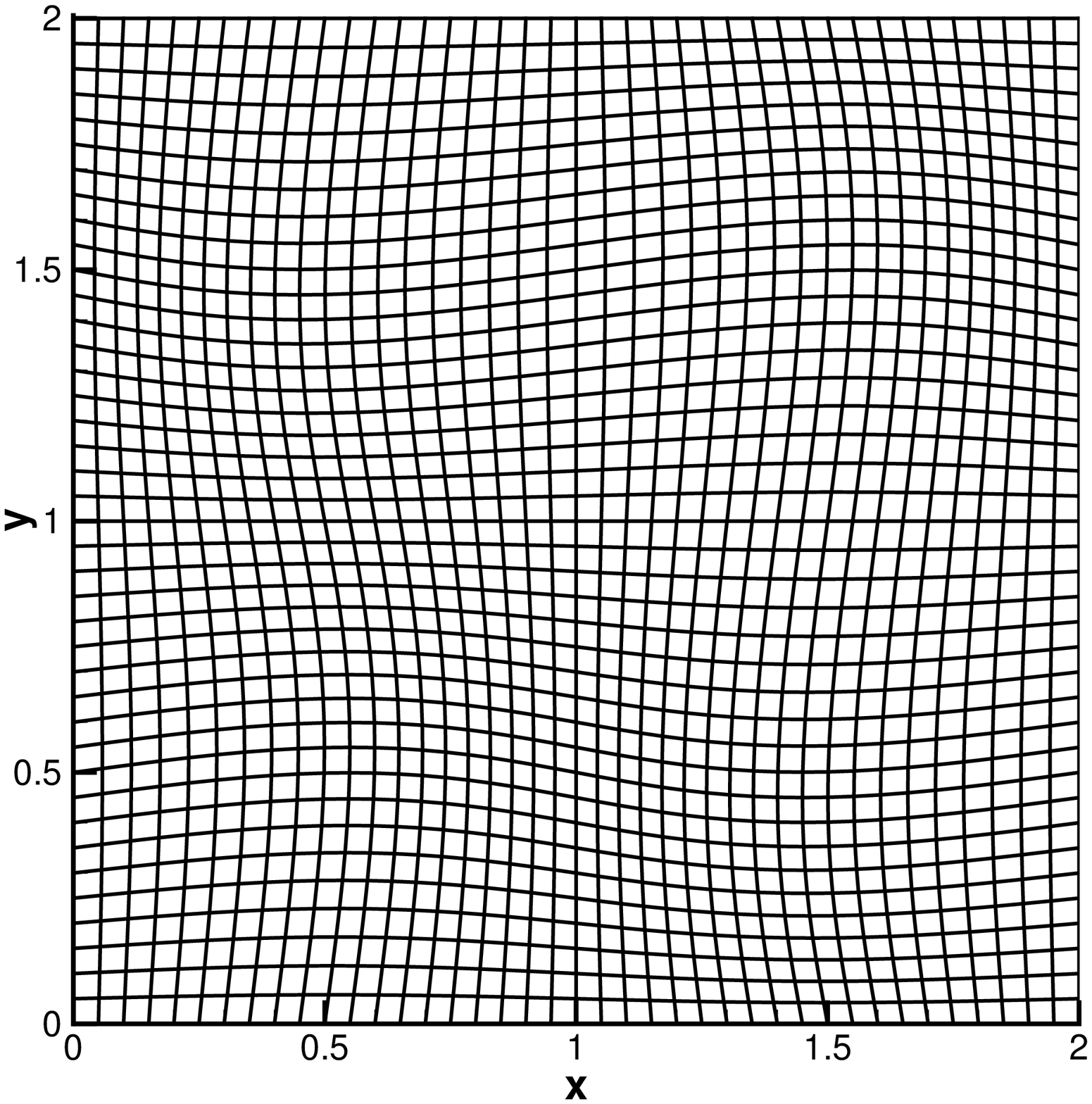}{b}
\caption{\label{accu-mesh} Accuracy test: 2D nonuniform orthogonal
mesh (a) and nonorthogonal mesh (b).}
\end{figure}

\begin{table}[!h]
\begin{center}
\def\temptablewidth{0.75\textwidth}
{\rule{\temptablewidth}{1.0pt}}
\begin{tabular*}{\temptablewidth}{@{\extracolsep{\fill}}c|cc|cc}
mesh      &  $L^1$ error &    order  &  $L^2$ error  &   order  \\
\hline
$10^2$    &  9.3427E-03  &   ~~      &   5.3475E-03  &  ~~      \\
$20^2$    &  3.0710E-04  &   4.9270  &   1.7197E-04  &  4.9585  \\
$40^2$    &  9.7453E-06  &   4.9778  &   5.4578E-06  &  4.9777  \\
$80^2$    &  3.0561E-07  &   4.9949  &   1.7120E-07  &  4.9945  \\
$160^2$   &  9.5621E-09  &   4.9982  &   5.3566E-09  &  4.9982  \\
$320^2$   &  2.9939E-10  &   4.9972  &   1.6773E-10  &  4.9971  \\
\end{tabular*}
{\rule{\temptablewidth}{1.0pt}}
\end{center}
\vspace{-5mm} \caption{\label{tab-2d-1} Accuracy test:
2D  advection of density perturbation with nonuniform
orthogonal meshes.}
\begin{center}
\def\temptablewidth{0.75\textwidth}
{\rule{\temptablewidth}{1.0pt}}
\begin{tabular*}{\temptablewidth}{@{\extracolsep{\fill}}c|cc|cc}
mesh      &  $L^1$ error &    order    &  $L^2$ error  &    order  \\
\hline
$10^2$    &  2.0481E-02  &       ~~    &  1.2061E-02   &      ~~   \\
$20^2$    &  9.4574E-04  &   4.43672   &  5.3232E-04   &    4.5019 \\
$40^2$    &  3.1638E-05  &   4.90167   &  1.7690E-05   &    4.9112 \\
$80^2$    &  9.9865E-07  &   4.98556   &  5.6017E-07   &    4.9809 \\
$160^2$   &  3.1326E-08  &   4.99453   &  1.7701E-08   &    4.9839 \\
$320^2$   &  9.8696E-10  &   4.98823   &  5.7052E-10   &    4.9554 \\
\end{tabular*}
{\rule{\temptablewidth}{1.0pt}}
\end{center}
\vspace{-5mm} \caption{\label{tab-2d-2} Accuracy test:
2D advection of density perturbation with nonuniform
nonorthogonal meshes.}
\begin{center}
\def\temptablewidth{0.75\textwidth}
{\rule{\temptablewidth}{1.0pt}}
\begin{tabular*}{\temptablewidth}{@{\extracolsep{\fill}}c|cc|cc}
mesh     & $L^1$ error  &    order   &  $L^2$ error  &   order  \\
\hline
$10^2$   &  6.7911E-03  &    ~~      &   3.7248E-03  &    ~~    \\
$20^2$   &  2.2028E-04  &   4.9462   &   1.2245E-04  &   4.9268 \\
$40^2$   &  7.0197E-06  &   4.9718   &   3.8911E-06  &   4.9758 \\
$80^2$   &  2.2575E-07  &   4.9585   &   1.2518E-07  &   4.9580 \\
$160^2$  &  7.7220E-09  &   4.8696   &   4.2895E-09  &   4.8670 \\
$320^2$  &  3.0956E-10  &   4.6406   &   1.7222E-10  &   4.6384 \\
\end{tabular*}
{\rule{\temptablewidth}{1.0pt}}
\end{center}
\vspace{-5mm} \caption{\label{tab-2d-3} Accuracy test:
2D advection of density perturbation with uniform
meshes.}
\end{table}

For the two-dimensional case, the physical domain is
$[0,2]\times[0,2]$ and the initial conditions are given as follows
\begin{align*}
\rho_0(x,y)=&1+0.2\sin(\pi (x+y)),~p_0(x,y)=1,\\
&U_0(x,y)=1,~V_0(x,y)=1.
\end{align*}
The periodic boundary conditions are imposed at boundaries and the
exact solutions are
\begin{align*}
\rho(x,y,t)=&1+0.2\sin(\pi(x+y-t)),~p(x,y,t)=1,\\
&U(x,y,t)=1,~V(x,y,t)=1.
\end{align*}
The computational domain is $[0,2]\times[0,2]$, and $N\times N$
uniform cells are used. The nonuniform orthogonal mesh and
nonuniform nonorthogonal mesh are tested respectively, where the
orthogonal is given by
\begin{align*}
\begin{cases}
\displaystyle x=\xi+0.05\sin (\pi \xi),\\
\displaystyle y=\eta+0.05\sin (\pi \eta),
\end{cases}
\end{align*}
the nonorthogonal mesh is given by
\begin{align*}
\begin{cases}
\displaystyle x=\xi+0.05\sin (\pi \xi)\sin (\pi \eta),\\
\displaystyle y=\eta+0.05\sin (\pi \xi)\sin (\pi \eta).
\end{cases}
\end{align*}
and the orthogonal and nonorthogonal meshes with $40\times40$ cells
are presented in Fig.\ref{accu-mesh} as example. As reference, the
uniform mesh with $N^2$ cells in the physical domain is also tested.
Two-dimensional Gauss quadratures are used to provide the initial
conditions. The $L^1$ and $L^2$ errors and orders of accuracy at
$t=2$ with $N^2$ cells are presented in Tab.\ref{tab-2d-1},
Tab.\ref{tab-2d-2} and Tab.\ref{tab-2d-3} for nonuniform orthogonal
meshes, nonorthogonal meshes and uniform meshes. The expected
accuracy can be also achieved for the current scheme.

\begin{table}[!h]
\begin{center}
\def\temptablewidth{0.75\textwidth}
{\rule{\temptablewidth}{1.0pt}}
\begin{tabular*}{\temptablewidth}{@{\extracolsep{\fill}}c|cc|cc}
mesh     & $L^1$ error  &    Order    &  $L^2$ error &  Order   \\
\hline
$10^3$   &  2.6560E-02  &    ~~       &  1.0639E-02   &   ~~    \\
$20^3$   &  9.0703E-04  &    4.8719   &  3.5650E-04   &  4.8993 \\
$40^3$   &  2.9298E-05  &    4.9522   &  1.1508E-05   &  4.9531 \\
$80^3$   &  9.5178E-07  &    4.9440   &  3.7407E-07   &  4.9432 \\
$160^3$  &  3.3343E-08  &    4.8351   &  1.3126E-08   &  4.8328 \\
\end{tabular*}
{\rule{\temptablewidth}{1.0pt}}
\end{center}
\vspace{-5mm}\caption{\label{tab-3d-1} Accuracy test: 3D advection
of density perturbation with nonuniform orthogonal meshes}
\begin{center}
\def\temptablewidth{0.75\textwidth}
{\rule{\temptablewidth}{1.0pt}}
\begin{tabular*}{\temptablewidth}{@{\extracolsep{\fill}}c|cc|cc}
mesh     & $L^1$ error  &   Order   &  $L^2$ error &  Order \\
\hline
$10^3$   &   3.5692E-02 &     ~~    &  1.4482E-02 &  ~~     \\
$20^3$   &   1.4233E-03 &   4.6482  &  6.1955E-04 &  4.5469 \\
$40^3$   &   4.9497E-05 &   4.8458  &  2.1579E-05 &  4.8435 \\
$80^3$   &   1.6323E-06 &   4.9223  &  7.0200E-07 &  4.9420 \\
$160^3$  &   5.6847E-08 &   4.8437  &  2.3694E-08 &  4.8888 \\
\end{tabular*}
{\rule{\temptablewidth}{1.0pt}}
\end{center}
\vspace{-5mm}\caption{\label{tab-3d-2} Accuracy test: 3D advection
of density perturbation with nonuniform nonorthogonal meshes}
\begin{center}
\def\temptablewidth{0.75\textwidth}
{\rule{\temptablewidth}{1.0pt}}
\begin{tabular*}{\temptablewidth}{@{\extracolsep{\fill}}c|cc|cc}
  mesh   & $L^1$ error  &     Order   &  $L^2$ error &  Order    \\
\hline
$10^3$   &   1.8756E-02  &  ~~      &     7.4373E-03  &      ~~  \\
$20^3$   &   6.3946E-04  &  4.8743  &     2.5037E-04  &   4.8926 \\
$40^3$   &   2.0744E-05  &  4.9460  &     8.1505E-06  &   4.9410 \\
$80^3$   &   7.1379E-07  &  4.8610  &     2.8072E-07  &   4.8596 \\
$160^3$  &   2.8947E-08  &  4.6240  &     1.1391E-08  &   4.6231 \\
\end{tabular*}
{\rule{\temptablewidth}{1.0pt}}
\end{center}
\vspace{-5mm}\caption{\label{tab-3d-3} Accuracy test: 3D advection
of density perturbation with uniform meshes.}
\end{table}

The three-dimensional accuracy test is presented as well, which is
the start point of the simulation of complex flows with complicated
geometry. The physical domain is $[0,2]\times[0,2]\times[0,2]$ and
the initial condition is set as follows
\begin{align*}
\rho_0&(x, y, z)=1+0.2\sin(\pi(x+y+z)),~p_0(x,y,z)=1,\\
&U_0(x,y,z)=1,~V_0(x,y,z)=1,~W_0(x,y,z)=1.
\end{align*}
The periodic boundary conditions are applied at boundaries, and the
exact solutions are
\begin{align*}
\rho(x,y&,z,t)=1+0.2\sin(\pi(x+y+z-t)),~p(x,y,z,t)=1,\\
&U(x,y,z,t)=1,~V(x,y,z,t)=1,~W(x,y,z,t)=1.
\end{align*}
The computational domain is $[0,2]\times[0,2]\times[0,2]$. The
nonuniform orthogonal mesh and nonuniform nonorthogonal mesh are
tested respectively, where the orthogonal is given by
\begin{align*}
\begin{cases}
\displaystyle x=\xi+0.05\sin (\pi \xi),\\
\displaystyle y=\eta+0.05\sin (\pi \eta),\\
\displaystyle z=\zeta+0.05\sin (\pi \zeta),
\end{cases}
\end{align*}
the nonorthogonal mesh is given by
\begin{align*}
\begin{cases}
\displaystyle x=\xi+0.05\sin (\pi \xi)\sin (\pi \eta)\sin (\pi \zeta),\\
\displaystyle y=\eta+0.05\sin (\pi \xi)\sin (\pi \eta)\sin (\pi \zeta),\\
\displaystyle z=\zeta+0.05\sin (\pi \xi)\sin (\pi \eta)\sin (\pi \zeta),
\end{cases}
\end{align*}
and $N^3$ uniform cells are used in the computational domain. As
reference, $N^3$ uniform cells in the physical domain is also
tested. The $L^1$ and $L^2$ errors and orders of accuracy at $t=2$
with $N^3$ cells are presented in Tab.\ref{tab-3d-1},
Tab.\ref{tab-3d-2} and Tab.\ref{tab-3d-3} for nonuniform orthogonal
meshes, nonorthogonal meshes, and uniform meshes. The expected
accuracy is also achieved for the current scheme for the
three-dimensional cases.

\begin{figure}[!htb]
\centering
\includegraphics[width=0.45\textwidth]{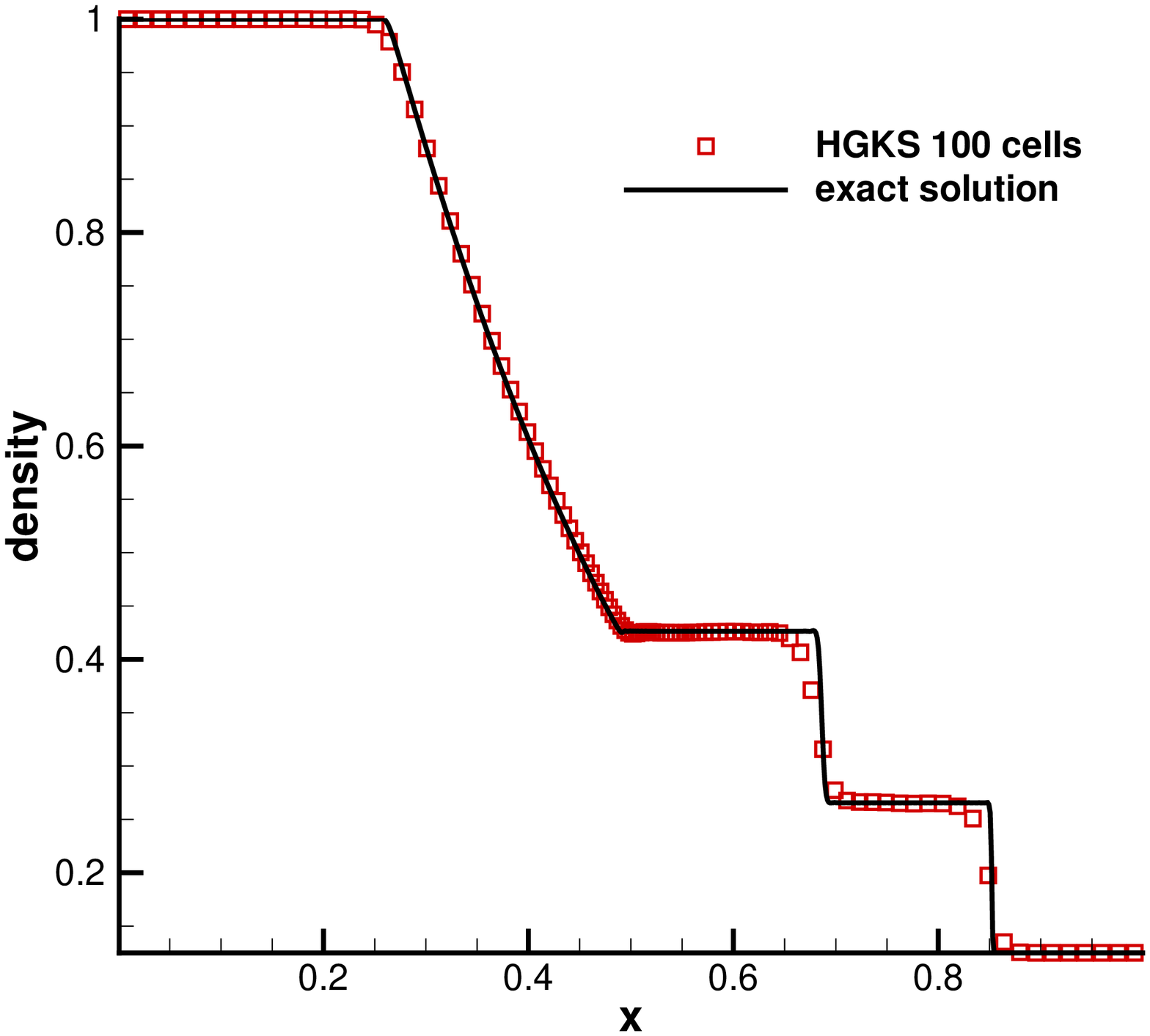}\includegraphics[width=0.45\textwidth]{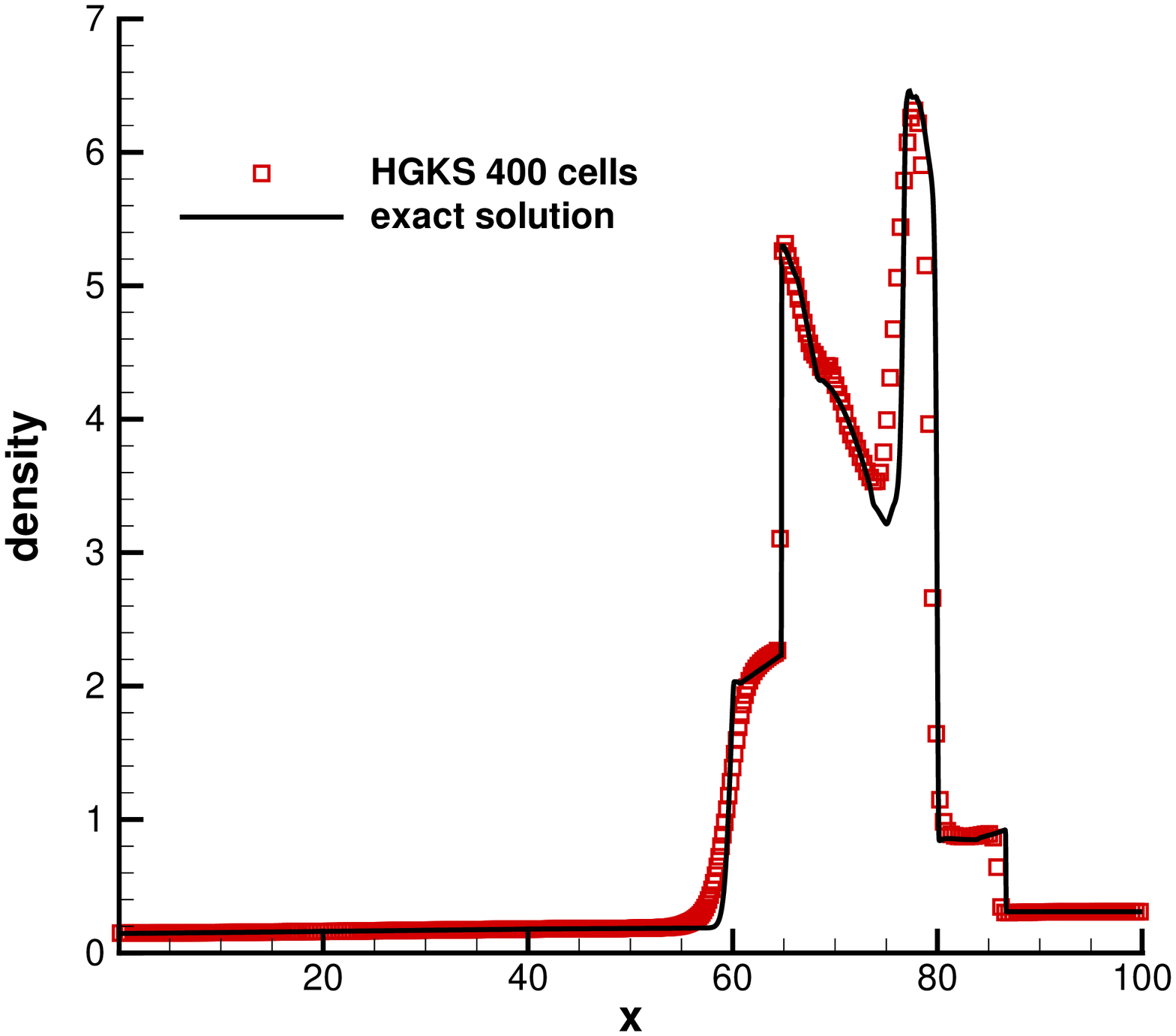}
\includegraphics[width=0.45\textwidth]{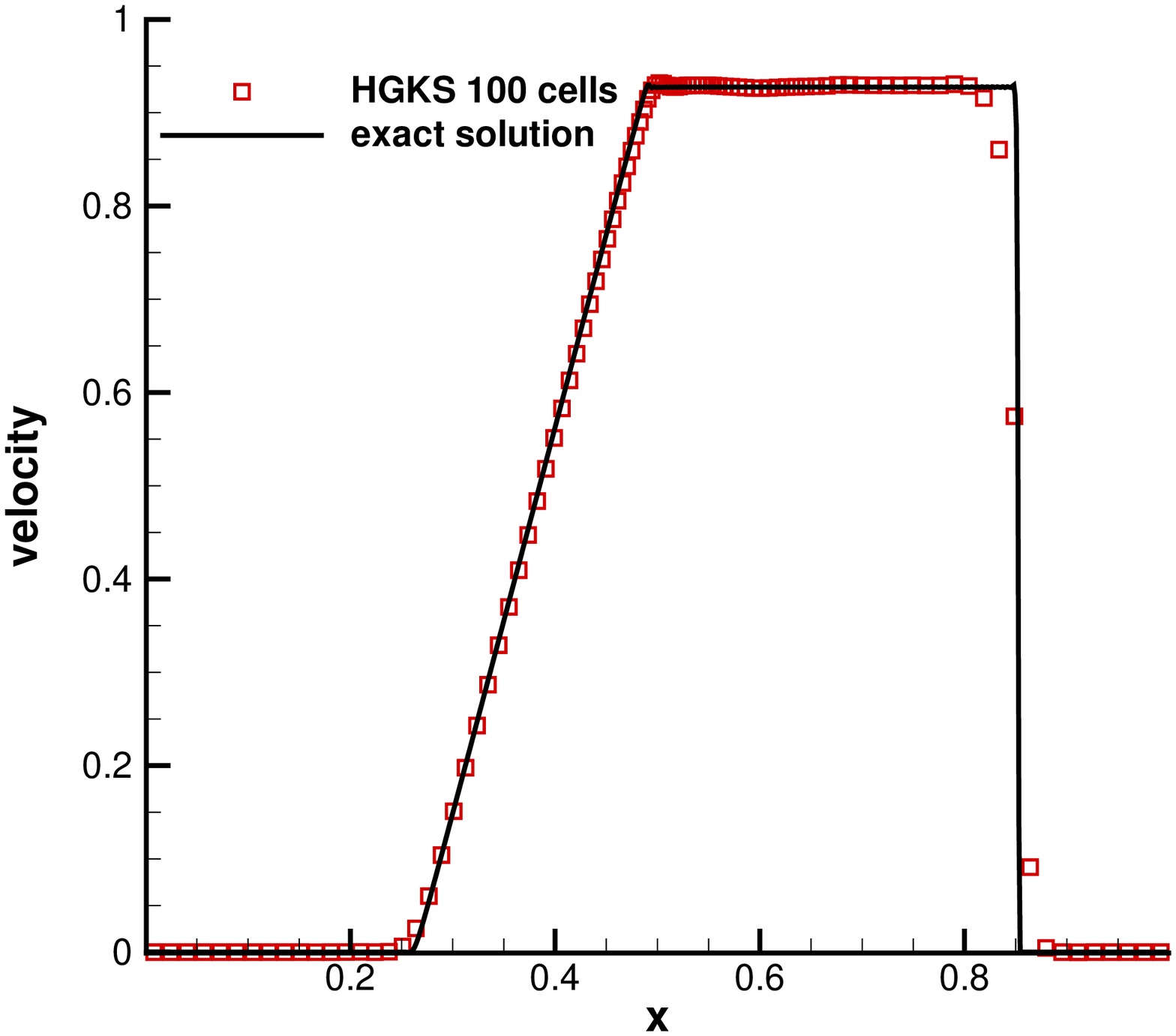}\includegraphics[width=0.45\textwidth]{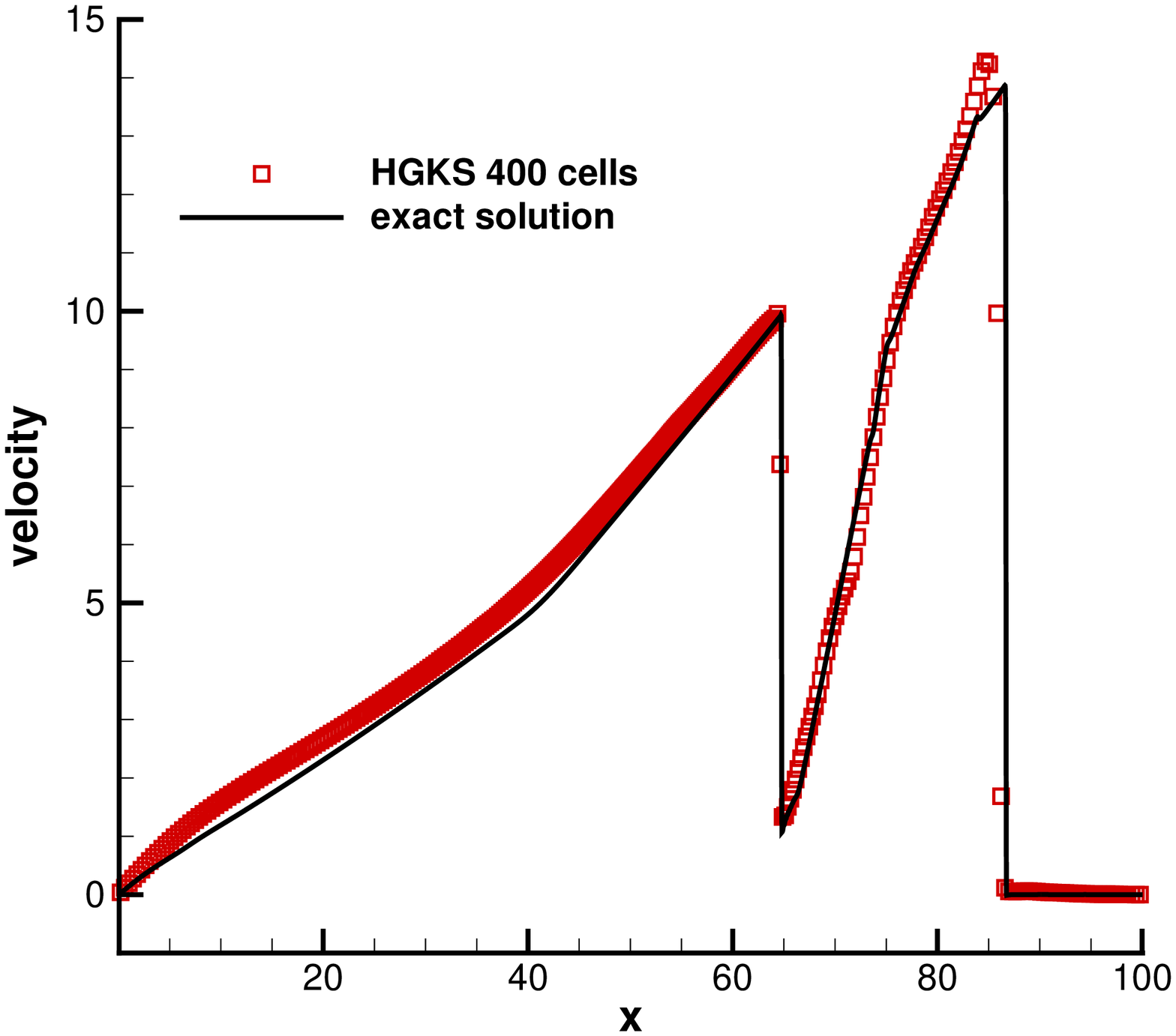}
\includegraphics[width=0.45\textwidth]{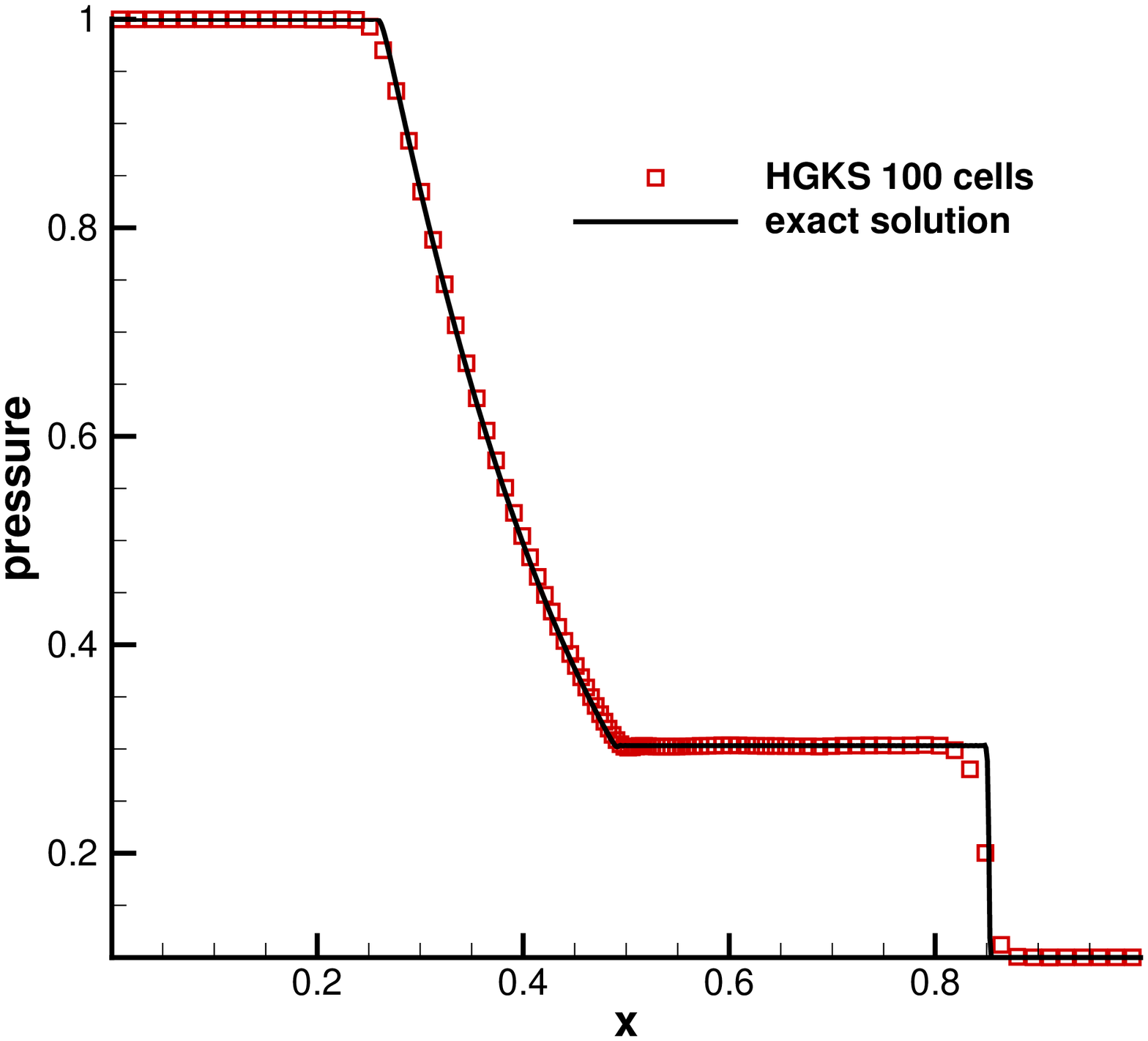}\includegraphics[width=0.45\textwidth]{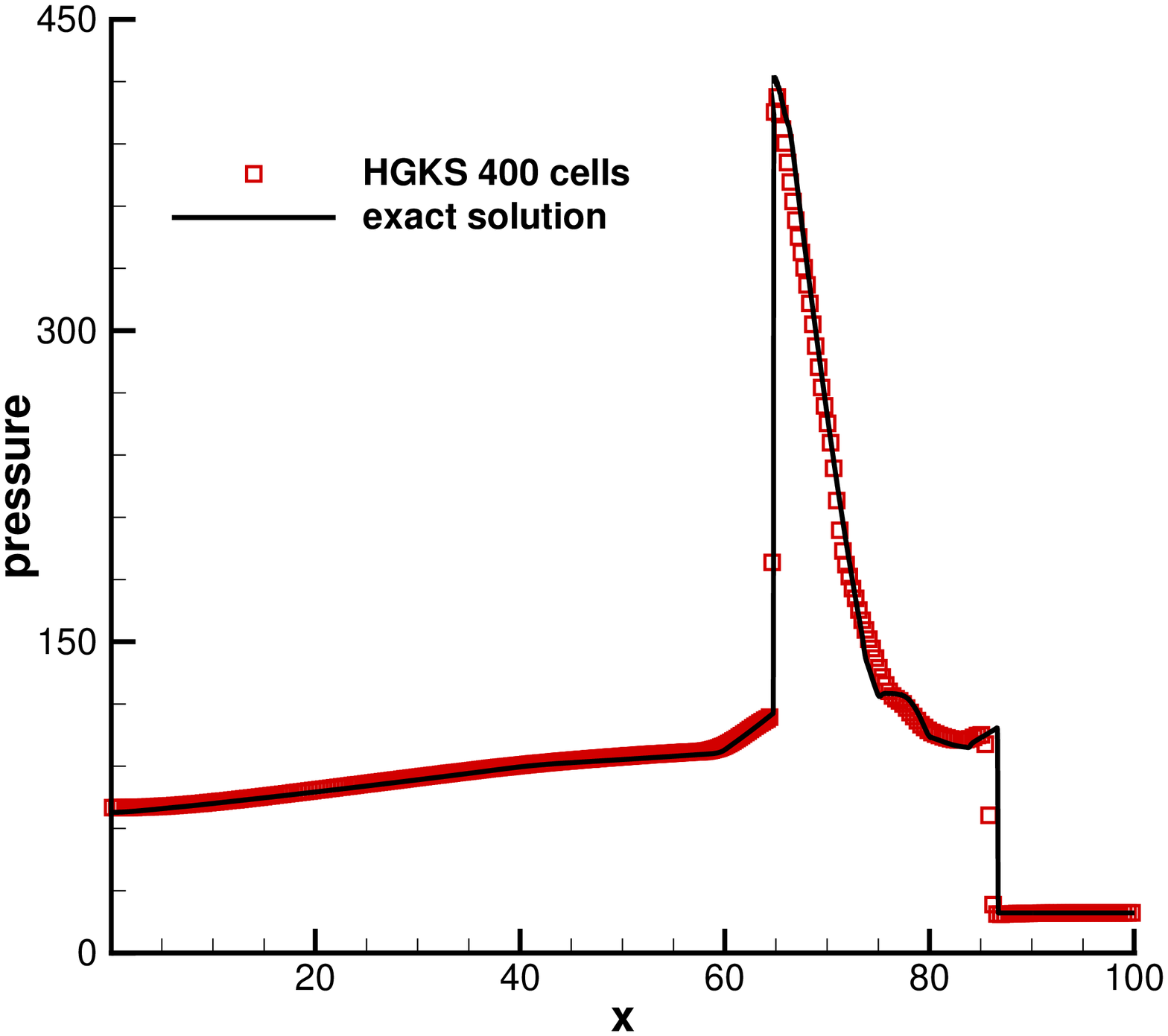}
\caption{\label{riemann-1}  One dimensional Riemann problem: the
density, velocity, pressure distributions for Sod problem at
$t=0.2$, and for blast wave problem at $t=0.038$.}
\end{figure}

\subsection{Geometric conservation law}
The geometric conservation law (GCL) \cite{GCL-1,GCL-2} is also
tested by the two-dimensional and three-dimensional nonuniform
nonorthogonal meshes given above. The GCL is mainly about the
maintenance of a uniform flow passing through a non-uniform
non-orthogonal mesh. The initial condition for the two-dimensional
case is
\begin{align*}
\rho_0(x,y)=1,~p_0(x,y)=1,~U_0(x,y)=1,~V_0(x,y)=1,
\end{align*}
and the initial condition for the three-dimensional case is
\begin{align*}
\rho_0&(x, y, z)=1,~p_0(x,y,z)=1,~U_0(x,y,z)=1,~V_0(x,y,z)=1,~W_0(x,y,z)=1.
\end{align*}
The periodic boundary conditions are adopted as well. The $L^1$ and
$L^2$ errors at $t=0.5$ for the two-dimensional case with $N^2$
cells are given in Tab.\ref{tab-2d-4}, and for the three-dimensional
case with $N^3$ cells are given in Tab.\ref{tab-3d-4}. The results
show that the errors reduce to the machine zero. The current scheme
is based on the coordinate transformation given by a smooth
function, which preserves the geometric conservation law
analytically. For a general mesh, the special treatment of the
metrics and Jacobian is needed \cite{curvilinear-6}.

\begin{table}[!h]
\begin{center}
\def\temptablewidth{0.5\textwidth}
{\rule{\temptablewidth}{1.0pt}}
\begin{tabular*}{\temptablewidth}{@{\extracolsep{\fill}}c|cc}
2D  mesh   & $L^1$ error  &  $L^2$ error   \\
\hline
$10^2$   &  2.9805E-15 &  2.3133E-15 \\
$20^2$   &  5.7204E-15 &  3.6866E-15 \\
$40^2$   &  7.5987E-15 &  4.7070E-15
\end{tabular*}
{\rule{\temptablewidth}{1.0pt}}
\end{center}
\vspace{-5mm} \caption{\label{tab-2d-4} Accuracy test:
two-dimensional geometric conservation law.}
\begin{center}
\def\temptablewidth{0.5\textwidth}
{\rule{\temptablewidth}{1.0pt}}
\begin{tabular*}{\temptablewidth}{@{\extracolsep{\fill}}c|cc}
3D   mesh   & $L^1$ error  &  $L^2$ error   \\
\hline
$10^3$   &  1.0896E-14 & 4.9111E-15\\
$20^3$   &  1.5292E-14 & 6.7811E-15\\
$40^3$   &  1.8241E-14 & 8.1087E-15
\end{tabular*}
{\rule{\temptablewidth}{1.0pt}}
\end{center}
\vspace{-5mm} \caption{\label{tab-3d-4} Accuracy test:
three-dimensional geometric conservation law.}
\end{table}

\subsection{One dimensional Riemann problem}
In this case, two examples of one-dimensional Riemann problems are
tested. The physical domain for the 1D case are $[0,1]$, and the
computational domain is expressed as
\begin{align*}
\displaystyle x=\xi+0.1\sin (2\pi \xi).
\end{align*}
The first one is the Sod problem, and the initial condition is given
as follows
\begin{equation*}
(\rho,U,p)=
\begin{cases}
(1,0,1),  0\leq x<0.5,\\
(0.125, 0, 0.1), 0.5\leq x\leq1.
\end{cases}
\end{equation*}
The non-reflecting boundary conditions are used at both ends, and
$100$ uniform cells are used in the computational domain. The
density, velocity and pressure distributions at $t=0.2$ are
presented in Fig.\ref{riemann-1}. The current scheme well captures
the exact solutions. The second one is the Woodward-Colella blast
wave problem \cite{Case-Woodward}, and the initial conditions are
given as follows
\begin{equation*}
(\rho,U,p)=\begin{cases}
(1, 0, 1000), \ \   0\leq x<0.1,\\
(1, 0, 0.01),\ \   0.1\leq x<0.9,\\
(1, 0, 100),\ \ \   0.9\leq x\leq 1.
\end{cases}
\end{equation*}
The reflected boundary conditions are imposed on both ends, and
$400$ non-uniform cells are used in the computational domain. The
density, velocity and pressure distributions at $t=0.038$ are
presented in Fig.\ref{riemann-1}, which validate the robustness and
resolution of currents scheme for the 1D strong discontinuity.

\begin{figure}[!htb]
\centering
\includegraphics[width=0.48\textwidth]{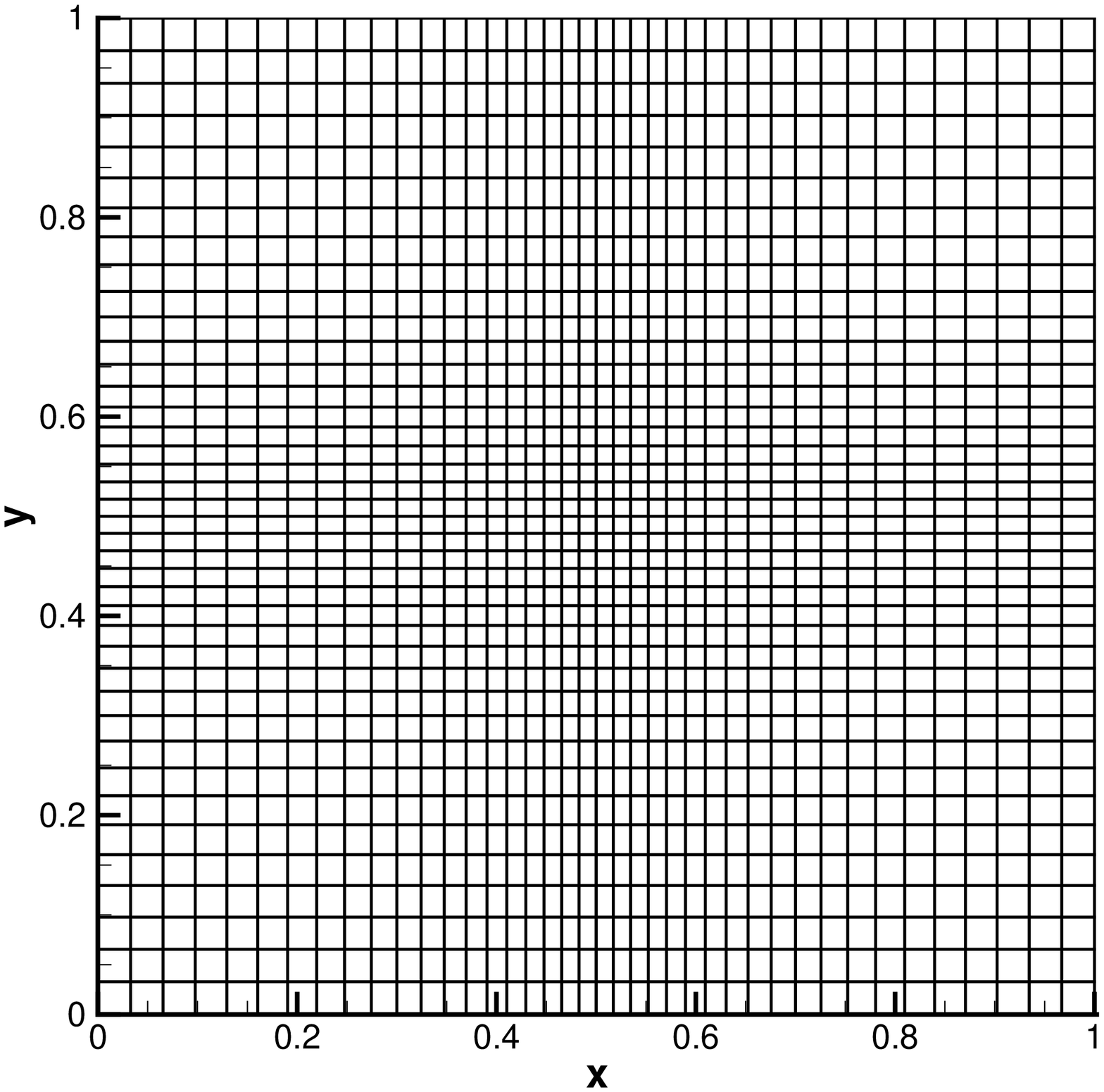}{a}
\includegraphics[width=0.48\textwidth]{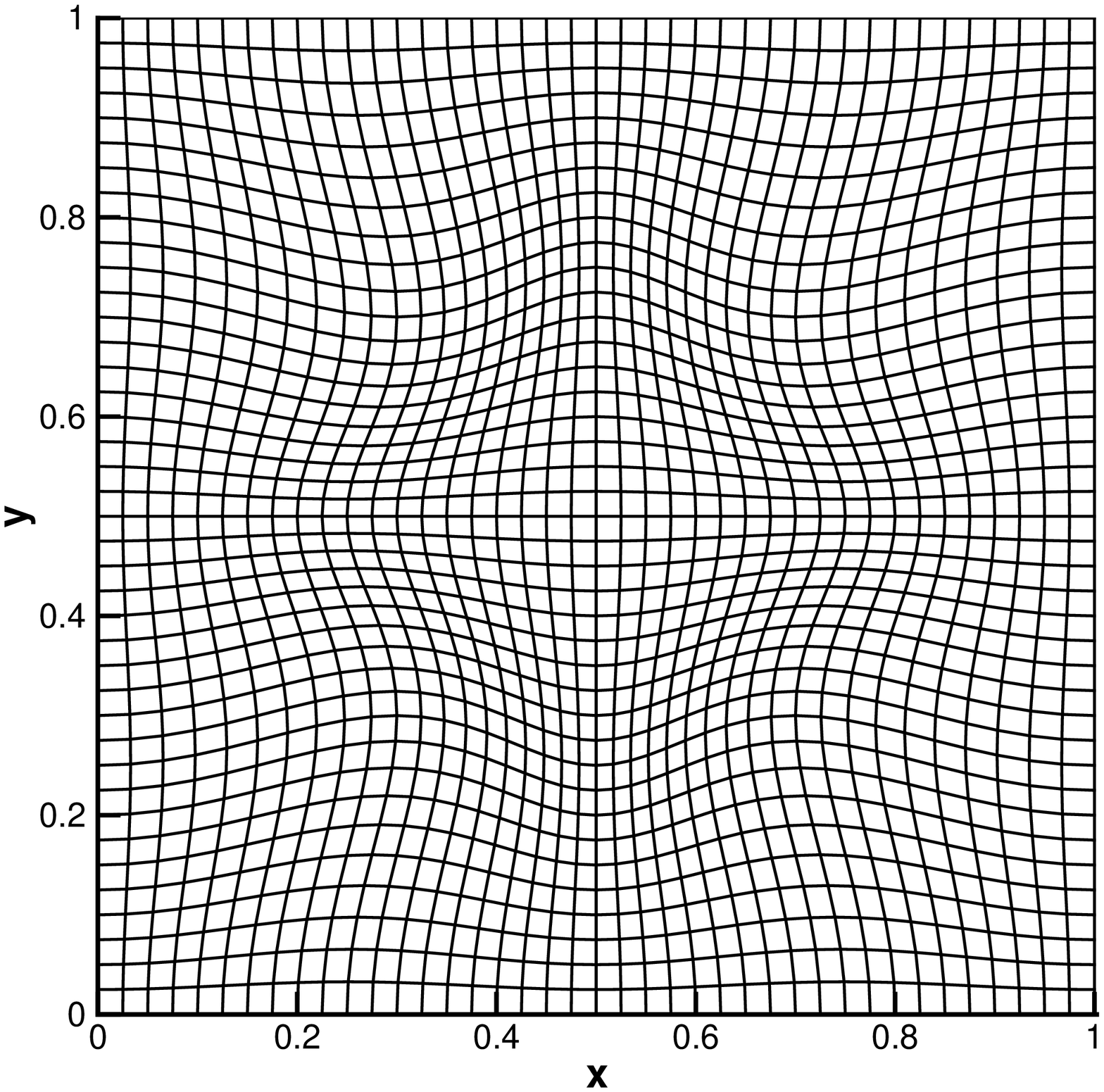}{b}
\caption{\label{Riemann-mesh}  2D Riemann problems: the nonuniform
orthogonal (a) and nonorthogonal (b) mesh $40\times40$ cells.}
\end{figure}

\begin{figure}[!htb]
\centering
\includegraphics[width=0.475\textwidth]{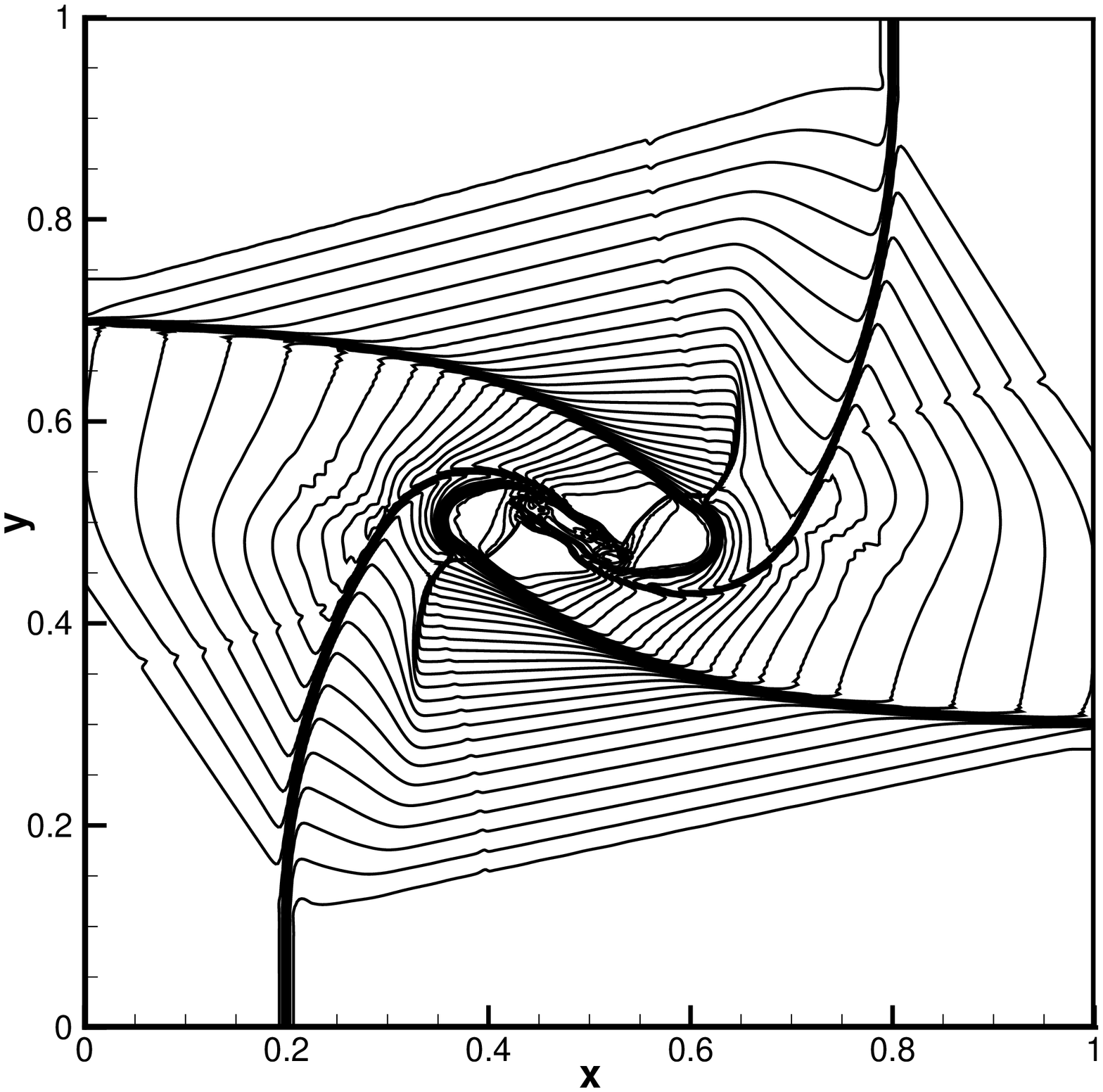}
\includegraphics[width=0.475\textwidth]{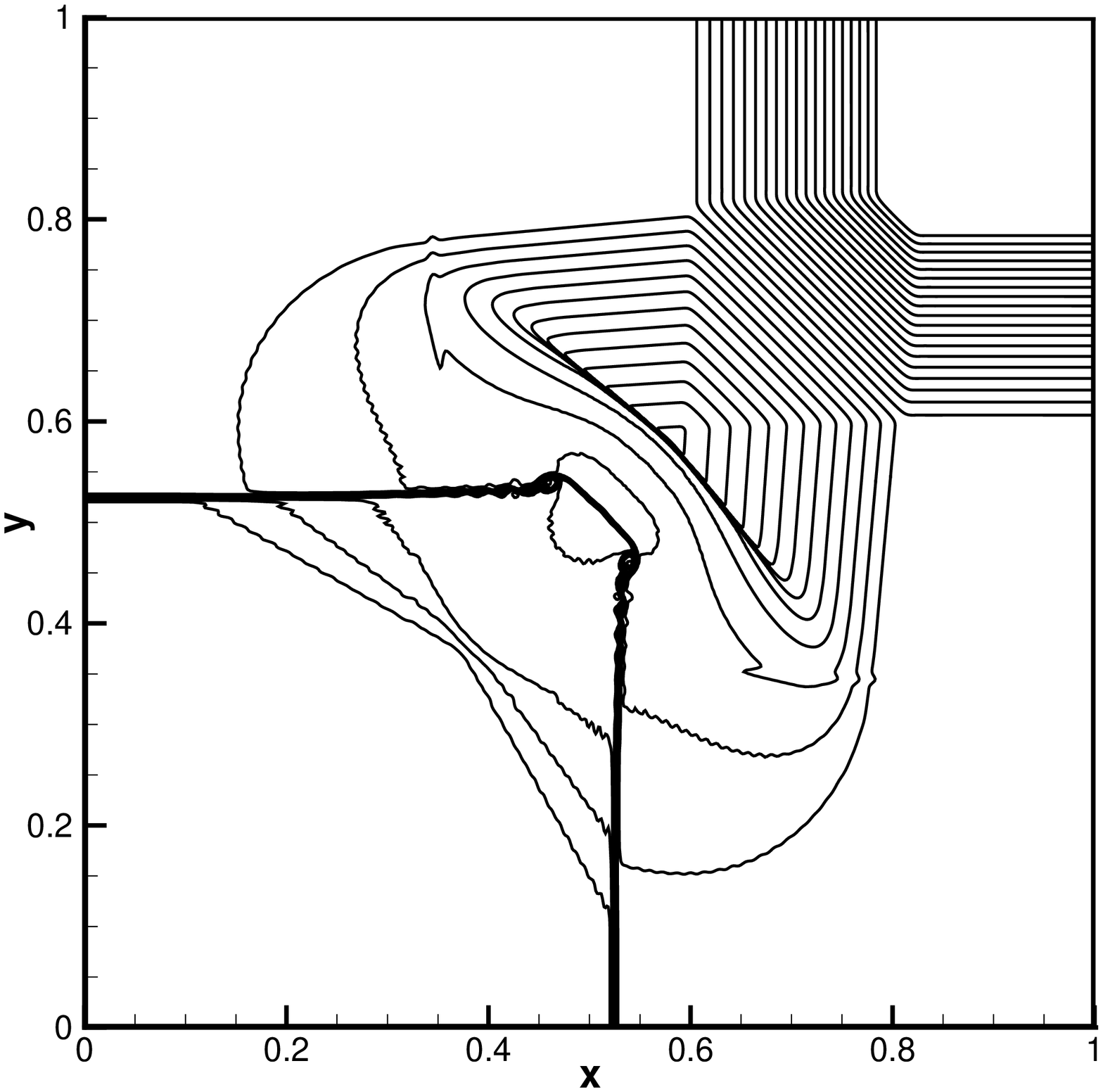}
\caption{\label{Riemann-2d-1} 2D Riemann problems: the density
distributions on the orthogonal mesh with $500\times500$ cells.}
\includegraphics[width=0.475\textwidth]{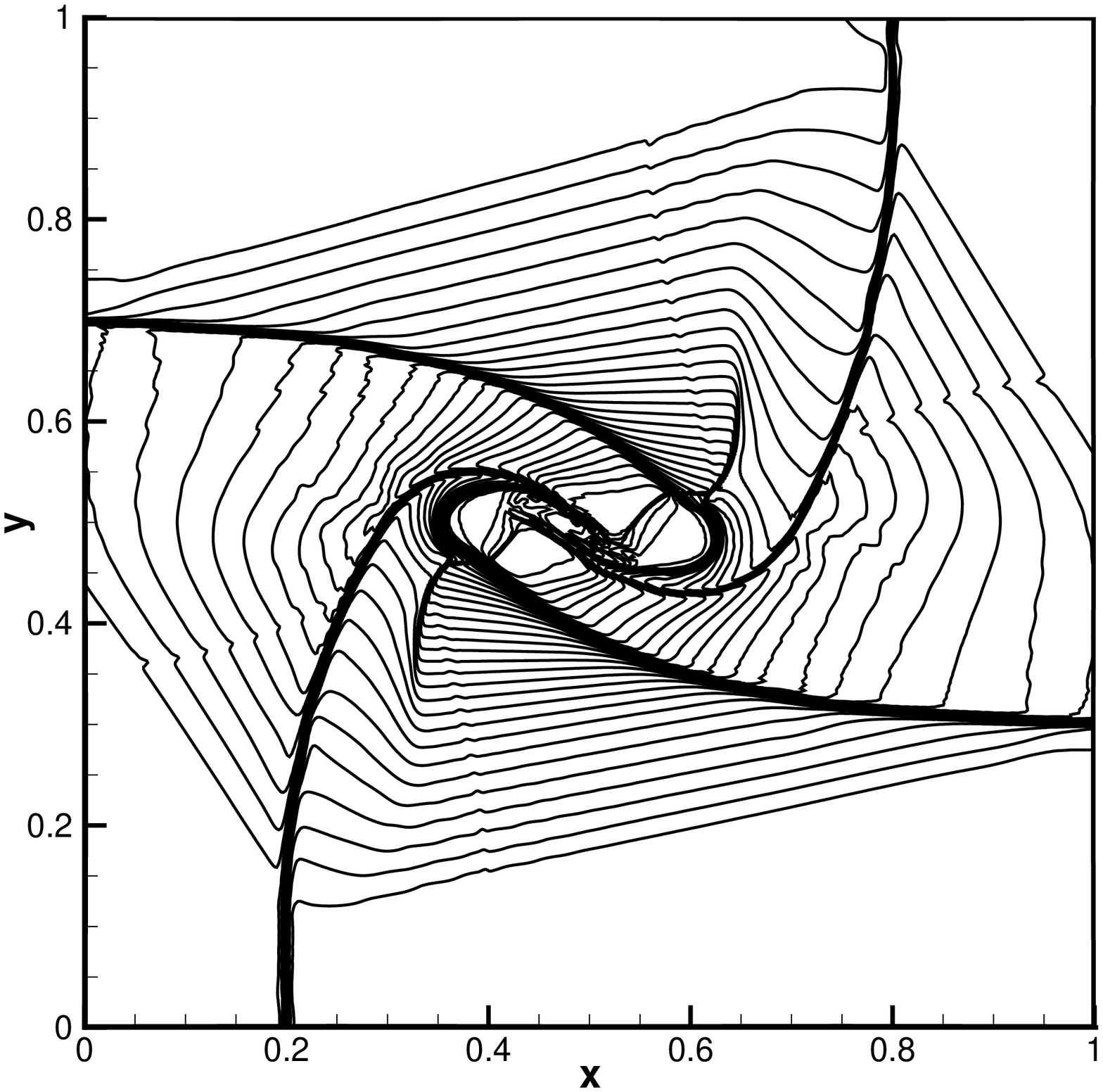}
\includegraphics[width=0.475\textwidth]{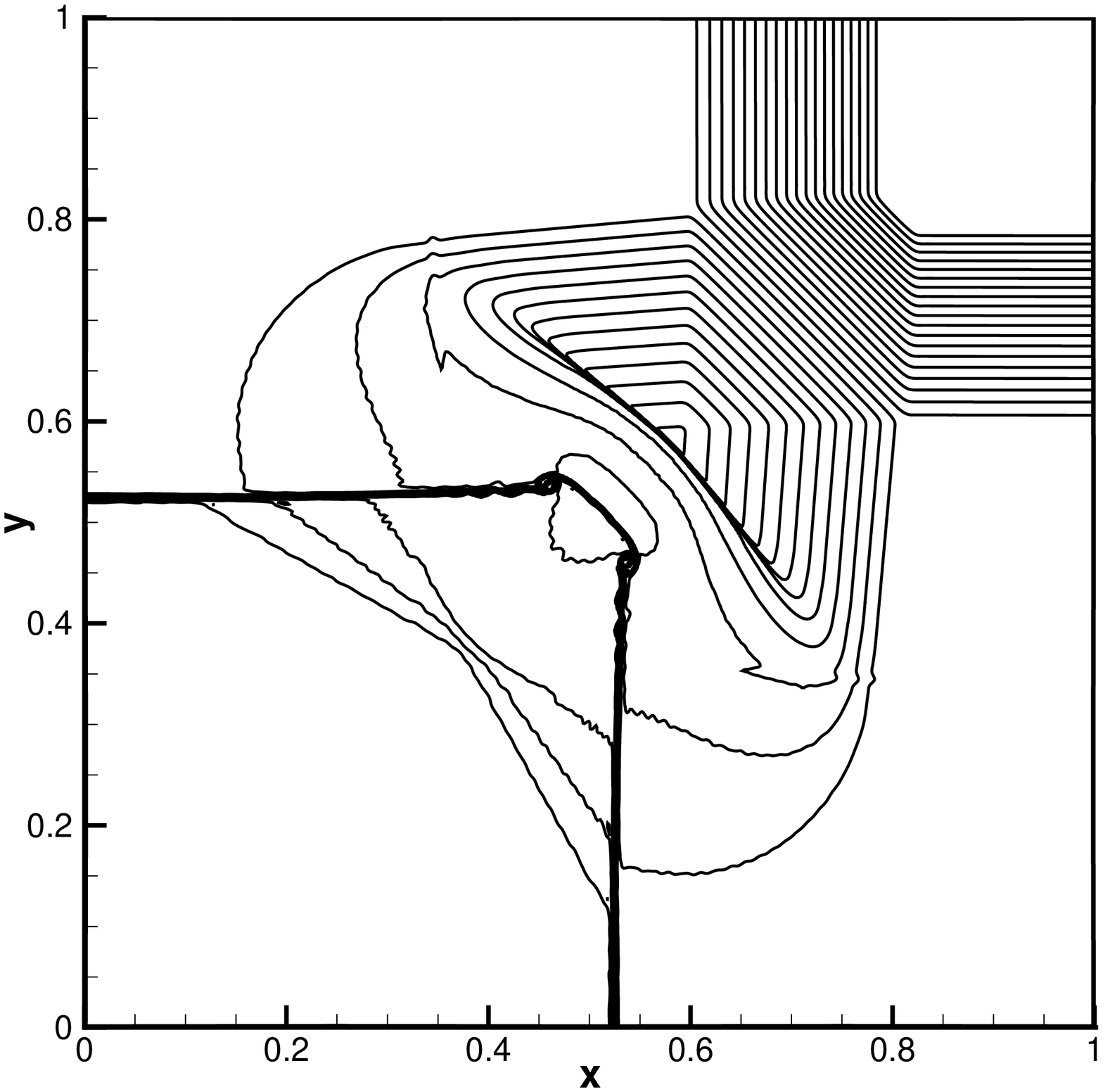}
\caption{\label{Riemann-2d-2} 2D Riemann problems: the density
distributions on the nonorthogonal mesh with $500\times500$ cells.}
\end{figure}

\subsection{Two-dimensional Riemann problems}
\begin{figure}[!htb]
\centering
\includegraphics[width=0.475\textwidth]{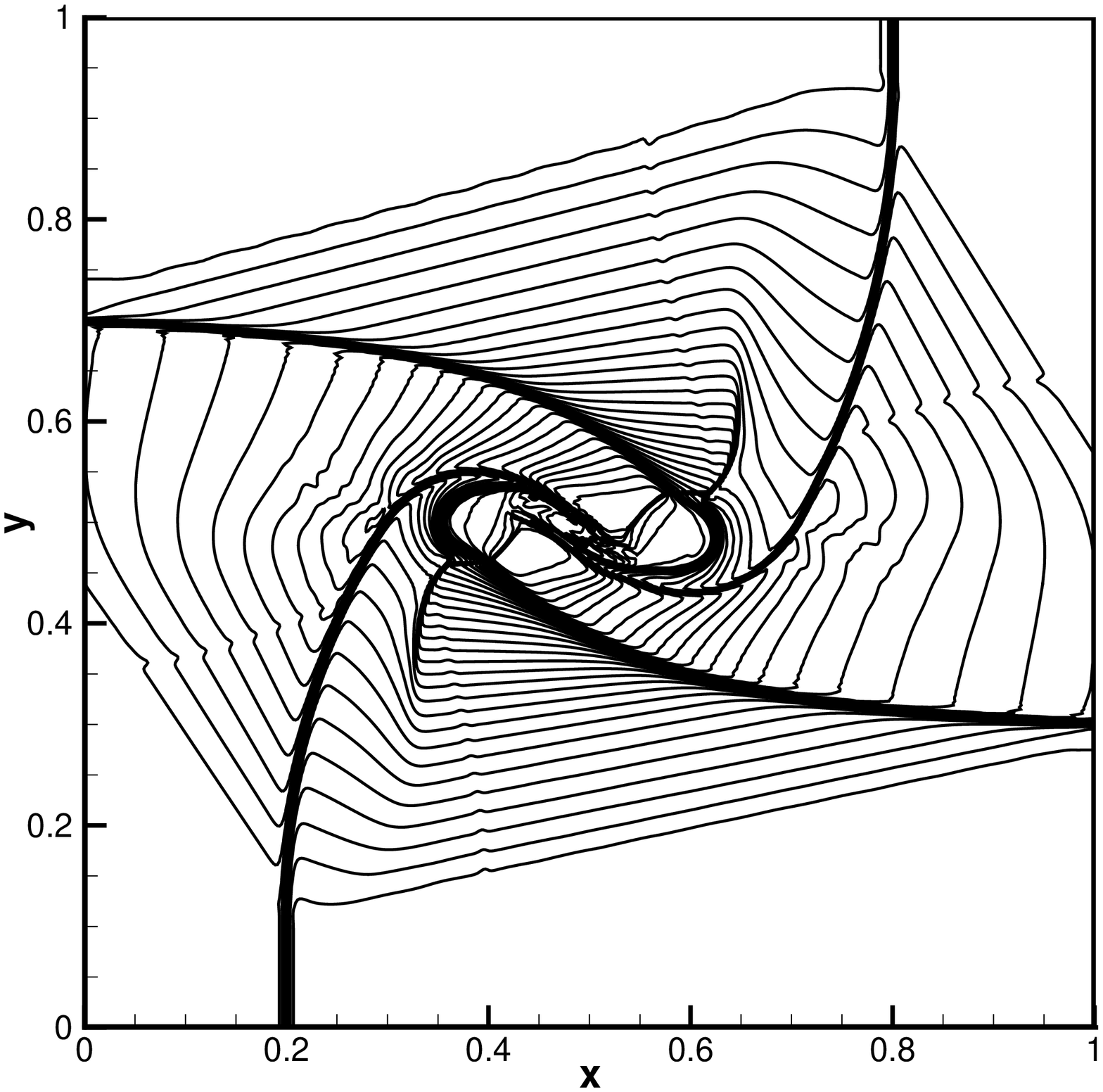}
\includegraphics[width=0.475\textwidth]{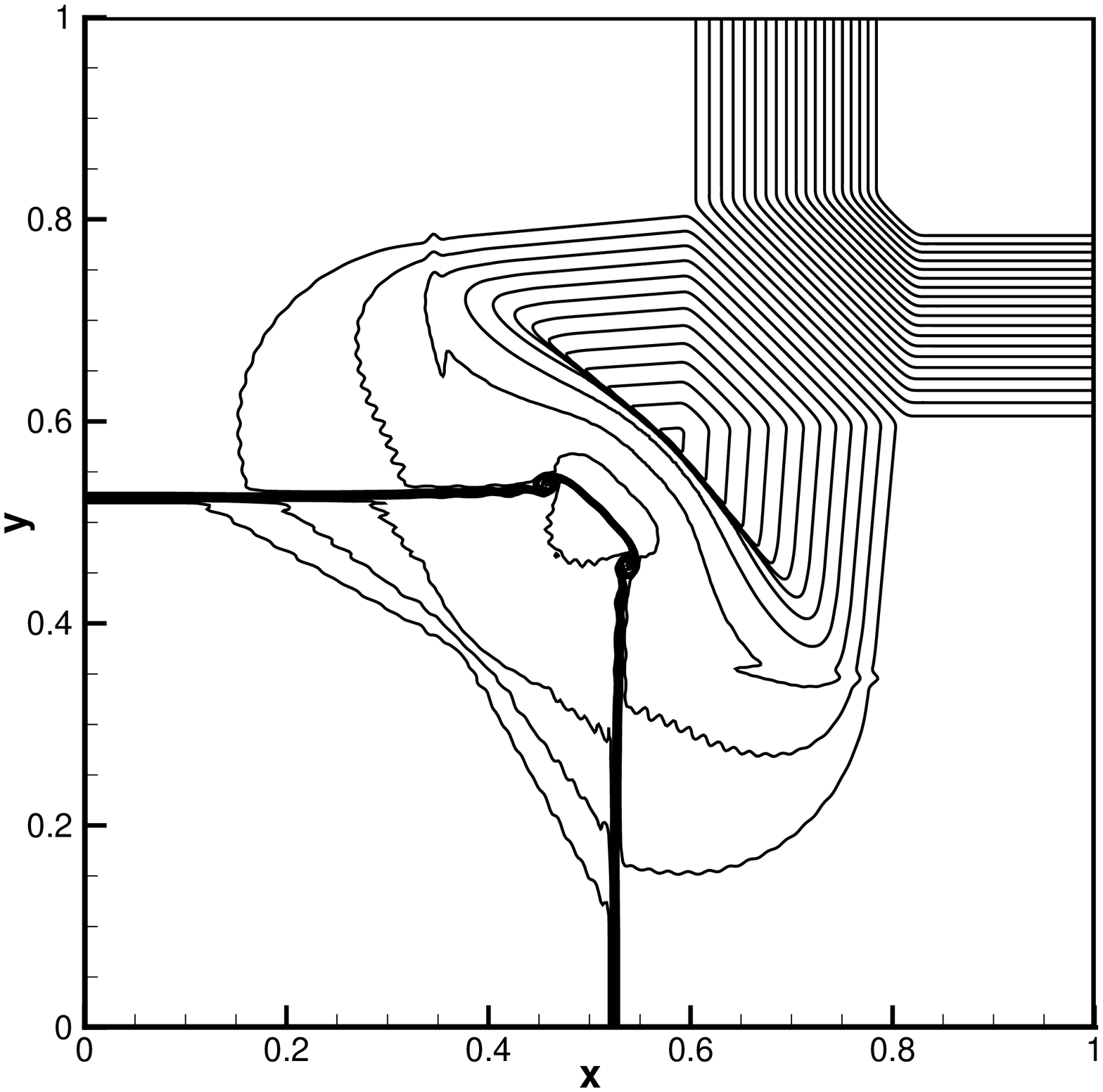}
\caption{\label{Riemann-2d-3} 2D Riemann problems: the density
distributions on the uniform mesh with $500\times500$ cells.}
\end{figure}

In this case, two examples of two-dimensional Riemann problems are
tested \cite{Case-Lax}. For these cases,  the nonuniform orthogonal
meshes
\begin{align*}
\begin{cases}
\displaystyle x=\xi+0.05\sin (2\pi \xi),\\
\displaystyle y=\eta+0.05\sin (2\pi \eta),
\end{cases}
\end{align*}
and nonuniform nonorthogonal meshes
\begin{align*}
\begin{cases}
\displaystyle x=\xi+0.05\sin^2 (2\pi \xi)\sin (2\pi \eta),\\
\displaystyle y=\eta+0.05\sin (2\pi \xi)\sin^2 (2\pi \eta).
\end{cases}
\end{align*}
are tested respectively, where both the physical and computational
domain are $[0,1]\times[0,1]$. The meshes with $40\times40$ cells
are shown in Fig.\ref{Riemann-mesh} as example. The initial
conditions for the first problem are
\begin{equation*}
(\rho,U,V,p) =\left\{\begin{array}{ll}
         (1,0.75,-0.5,0.5), \ \ \ &x>0.5,y>0.5,\\
         (2,0.75,0.5,0.5), &x<0.5,y>0.5,\\
         (1,-0.75£¬0.5,0.5), &x<0.5,y<0.5,\\
         (3,-0.75,-0.5, 0.5), &x>0.5,y<0.5,
                          \end{array} \right.
                          \end{equation*}
in which four initial contacts waves interact with each other and
result in a more complicated pattern. For the second case, the
initial conditions are
\begin{equation*}
(\rho,U,V,p)=\left\{\begin{aligned}
         &(1, 0.1, 0.1, 1), &x>0.5,y>0.5,\\
         &(0.5197,-0.6259, 0.1, 0.4), &x<0.5,y>0.5,\\
         &(0.8, 0.1, 0.1, 0.4), &x<0.5,y<0.5,\\
         &(0.5197,0.1,-0.6259, 0.4), &x>0.5,y<0.5,
                          \end{aligned} \right.
                          \end{equation*}
which simulate the interaction of the rarefaction waves and the
vortex-sheets. The non-reflecting boundary conditions are used in
all boundaries. Meanwhile, the meshes are given by symmetrically
corresponding the boundaries. The density distributions for the
first case at $t=0.4$ and for the second case at $t=0.25$ on the
orthogonal and nonorthogonal meshes with $500\times500$ cells are
presented in Fig.\ref{Riemann-2d-1} and Fig.\ref{Riemann-2d-2},
respectively. As reference, these two cases are tested on the
uniform mesh with $500\times500$ cells and the density distributions
are given in Fig.\ref{Riemann-2d-3}. The complicated flow structures
are well captures by the current scheme with different type of
meshes.

\begin{figure}[!htb]
\centering
\includegraphics[width=0.27\textwidth]{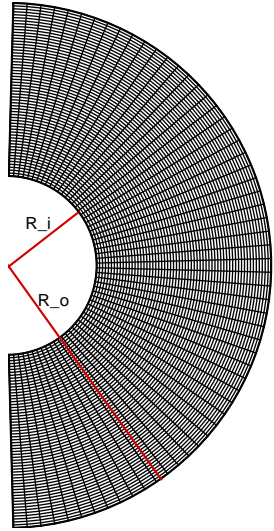}
\caption{\label{cylinder-1} Hypersonic flow past a cylinder: the
mesh with $60\times60$ cells.}
\end{figure}

\subsection{Hypersonic flow past a cylinder}
\begin{figure}[!htb]
\centering
\includegraphics[width=0.27\textwidth]{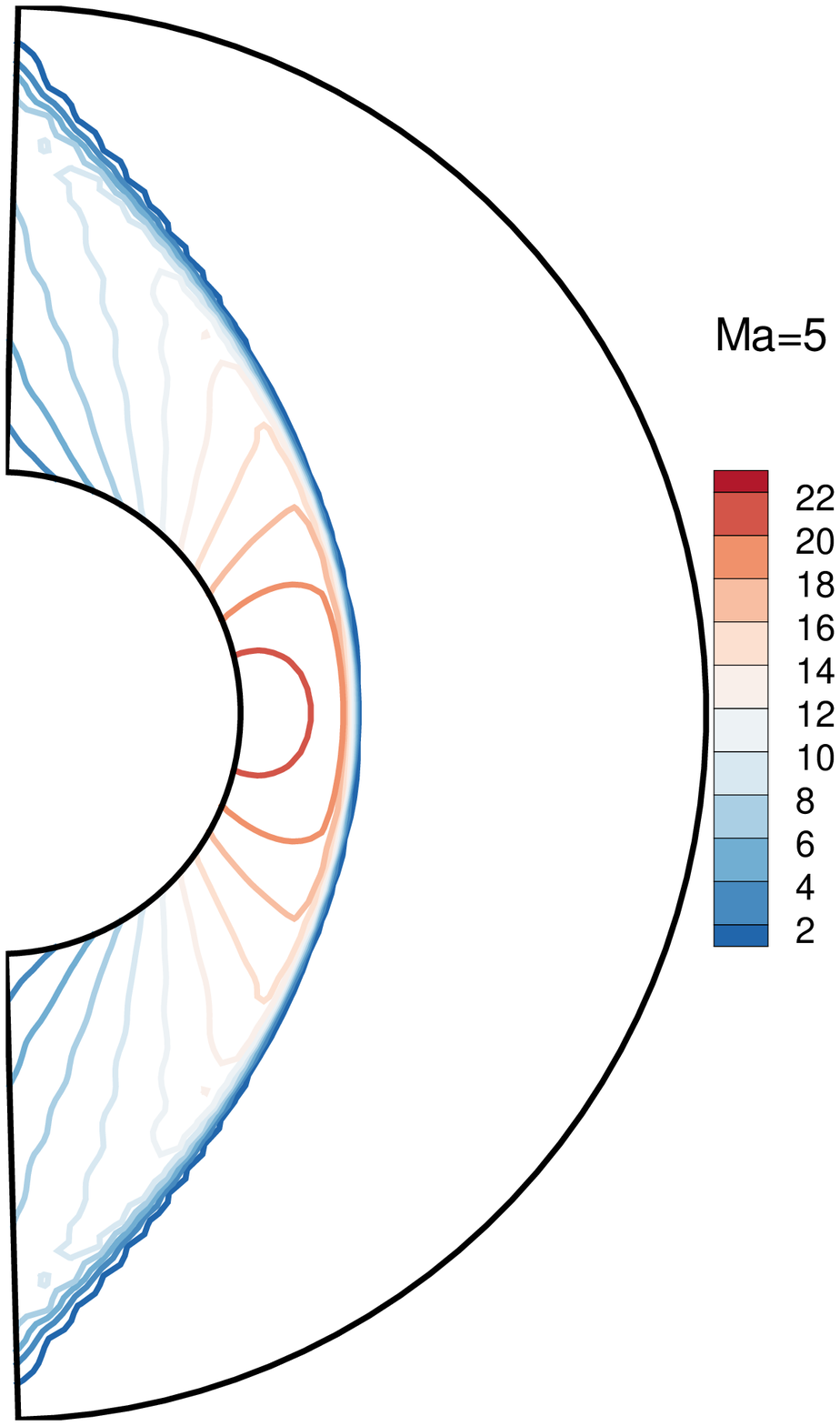}
\includegraphics[width=0.27\textwidth]{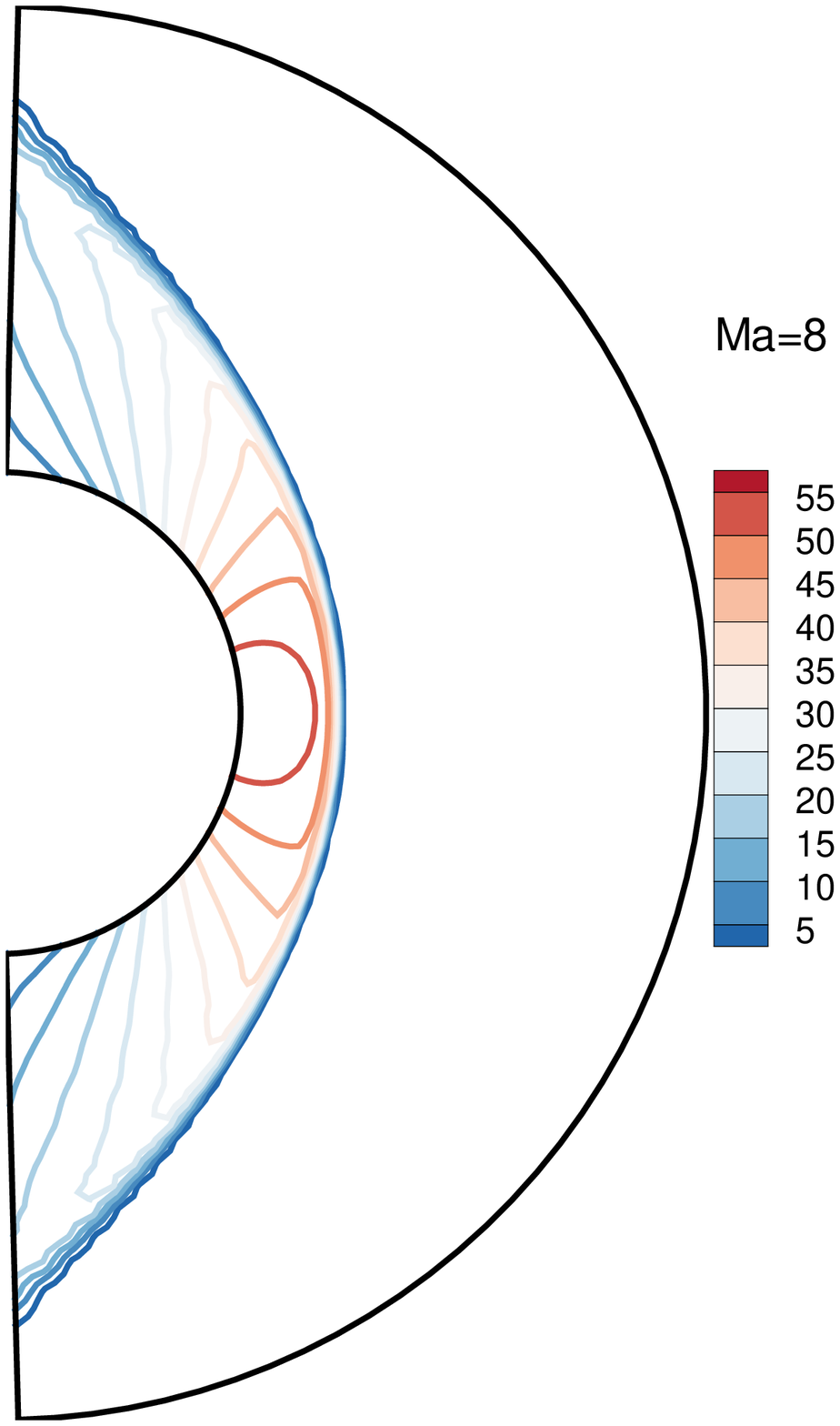}
\includegraphics[width=0.27\textwidth]{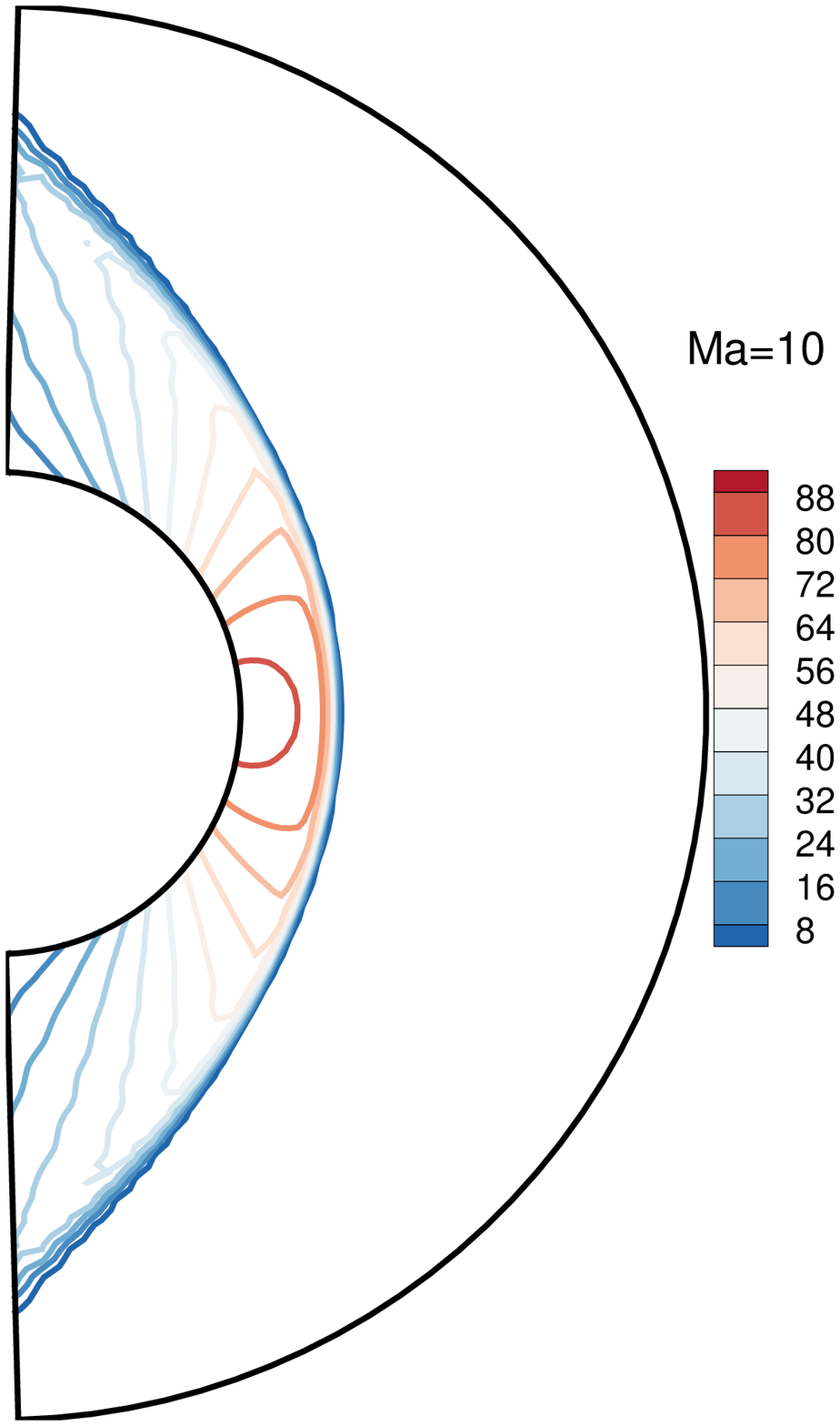}
\caption{\label{cylinder-2} Hypersonic flow past a cylinder: the
pressure distribution with Mach number $Ma=5, 8$ and $10$.}
\end{figure}

In this case, the hypersonic flows impinging on a cylinder are
tested to validate robustness of the current scheme for the inviscid
flow.  For this case, the computational domain is
$[0.5,1.5]\times[-0.5,0.5]$, and the physical domain is expressed as
\begin{align*}
\begin{cases}
\displaystyle x=\xi\cos(\pi\eta),\\
\displaystyle y=\xi\sin(\pi\eta),
\end{cases}
\end{align*}
In the computation, $60\times60$ cells are used shown in
Fig.\ref{cylinder-1}, which are given uniformly in the computational
domain. This problem is initialized by a flow moving towards to a
cylinder with different Mach numbers. The reflective boundary
condition is imposed on the surface of cylinder, and the outflow
boundary condition is set on the left boundary. The Mach number
distributions for the flows with $Ma=5, 8$, and $10$ are presented
in Fig.\ref{cylinder-2}, which show that the current scheme can
capture strong shocks very well without carbuncle phenomenon
\cite{Case-Pandolfi}. The robustness of the scheme is well
validated.

\subsection{Viscous shock tube}
This problem was introduced to test the performances of current
scheme for viscous flows \cite{Case-Daru}. In this case, an ideal
gas is at rest in a two-dimensional unit box $[0,1]\times[0,1]$. A
membrane located at $x=0.5$ separates two different states of the
gas and the dimensionless initial states are
\begin{equation*}
(\rho,U,p)=\left\{\begin{aligned}
&(120, 0, 120/\gamma), \ \ \ &  0<x<0.5,\\
&(1.2, 0, 1.2/\gamma),  & 0.5<x<1,
\end{aligned} \right.
\end{equation*}
where $\gamma=1.4$, Reynolds number $\mbox{Re} =200$ and Prandtl
number $\mbox{Pr} =0.73$. In the computation, this case is tested in
the physical domain $[0, 1]\times[0, 0.5]$, a symmetric boundary
condition is used on the top boundary $x\in[0, 1], y=0.5$.  Non-slip
boundary condition for velocity, and adiabatic condition for
temperature are imposed at solid wall boundaries. For this case, the
nonuniform orthogonal meshes
\begin{align*}
\begin{cases}
\displaystyle x=\xi-0.05\sin (2\pi \xi),\\
\displaystyle y=\eta-0.05\sin (2\pi \eta),
\end{cases}
\end{align*}
and nonuniform nonorthogonal meshes
\begin{align*}
\begin{cases}
\displaystyle x=\xi-0.05\sin^2 (2\pi \xi)\sin (2\pi \eta),\\
\displaystyle y=\eta-0.05\sin (2\pi \xi)\sin^2 (2\pi \eta).
\end{cases}
\end{align*}
are used, and the meshes with $50\times25$ cells are shown in Fig.\ref{Riemann-mesh} as example.

\begin{figure}[!h]
\centering
\includegraphics[width=0.55\textwidth]{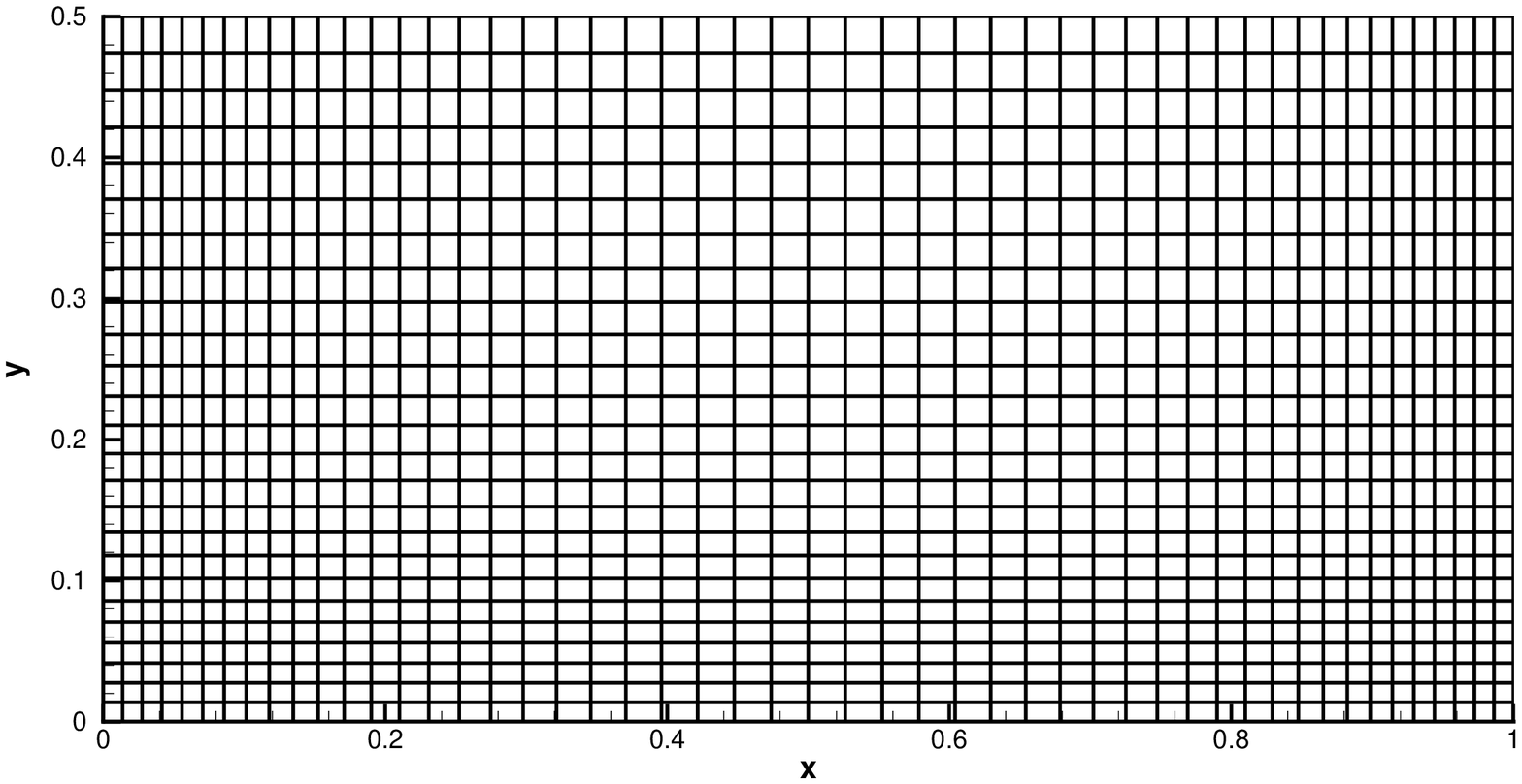}{a}
\includegraphics[width=0.55\textwidth]{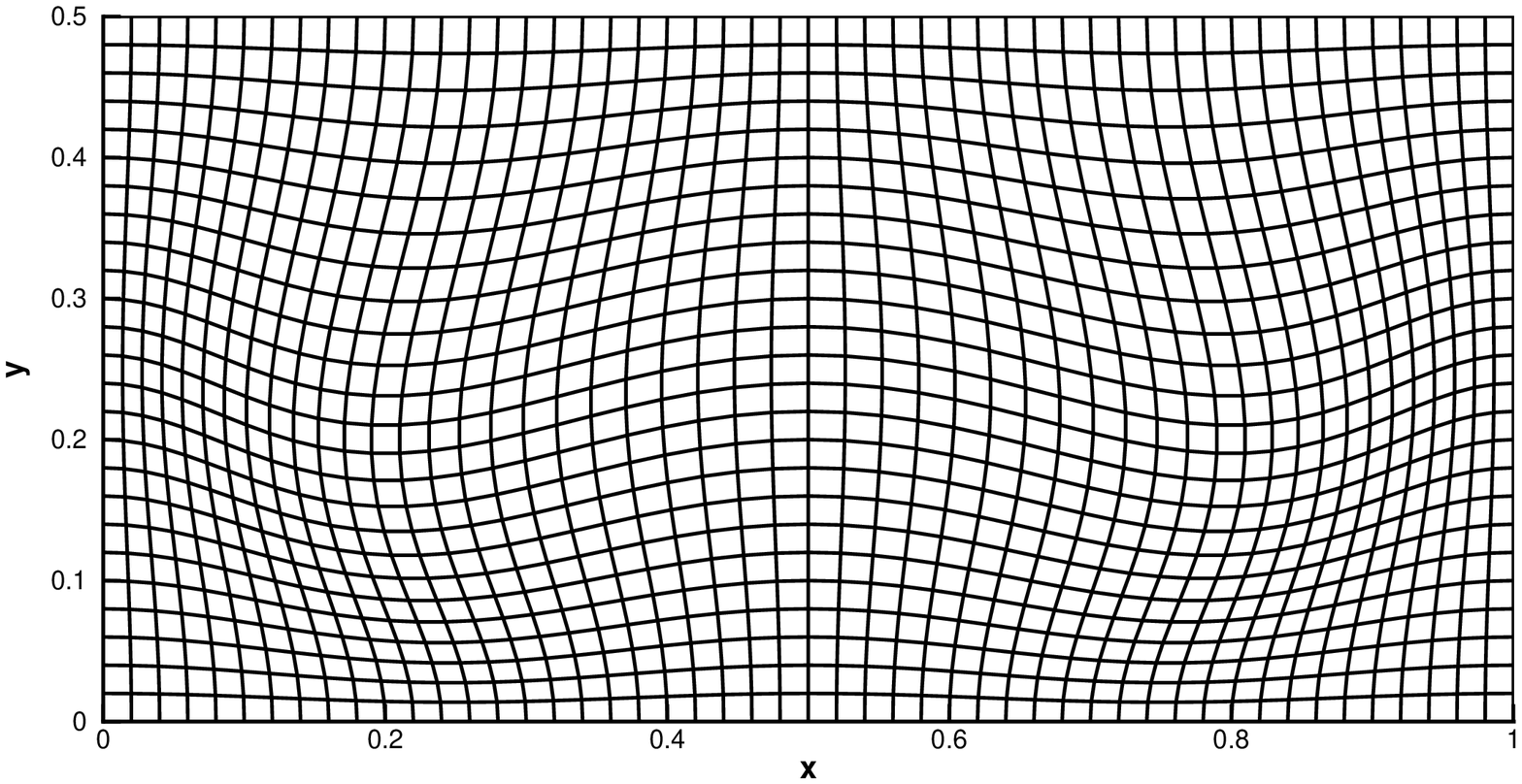}{b}
\caption{\label{vis-shock} Viscous shock tube: the nonuniform
orthogonal mesh (a) and nonorthogonal mesh (b).}
\end{figure}

\begin{figure}[!h]
\centering
\includegraphics[width=0.55\textwidth]{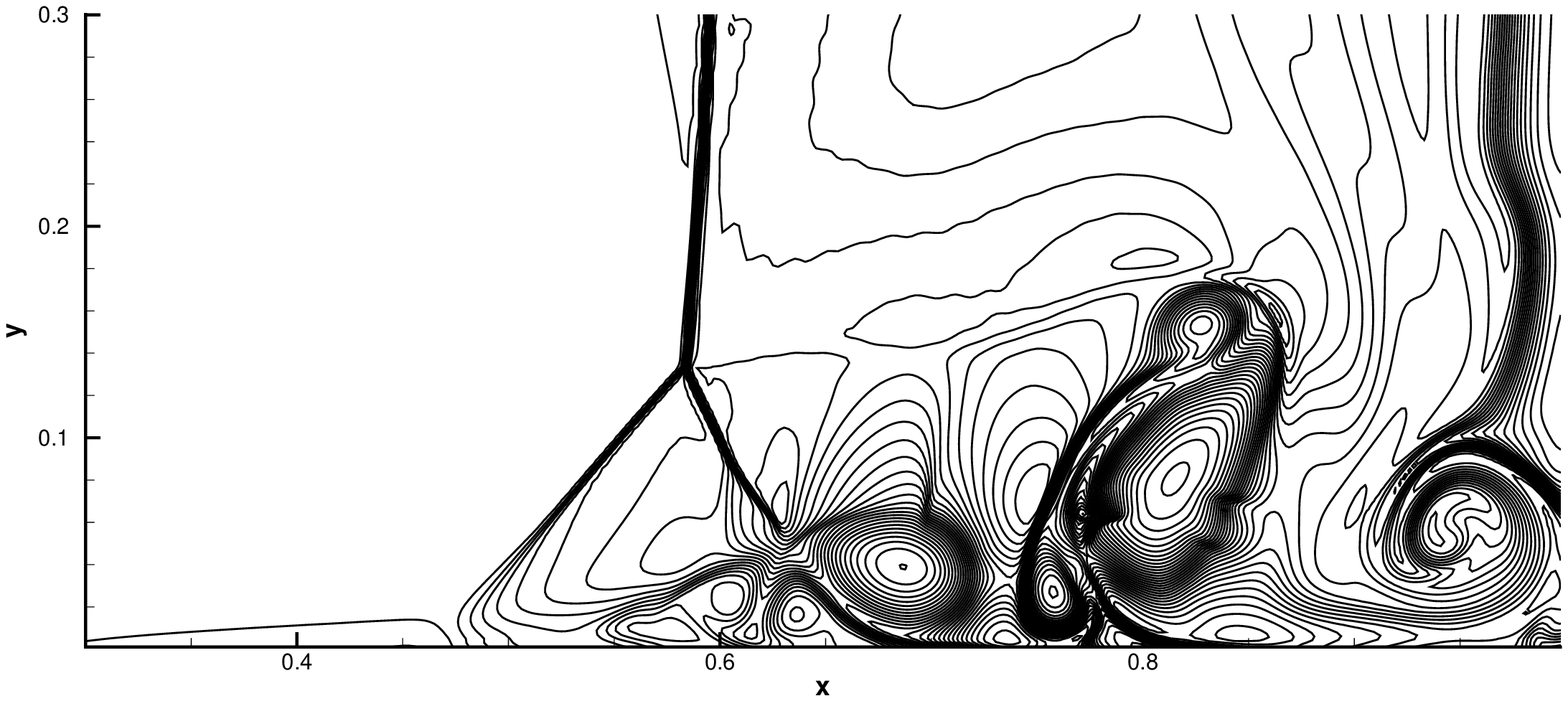}{a}\\
\includegraphics[width=0.55\textwidth]{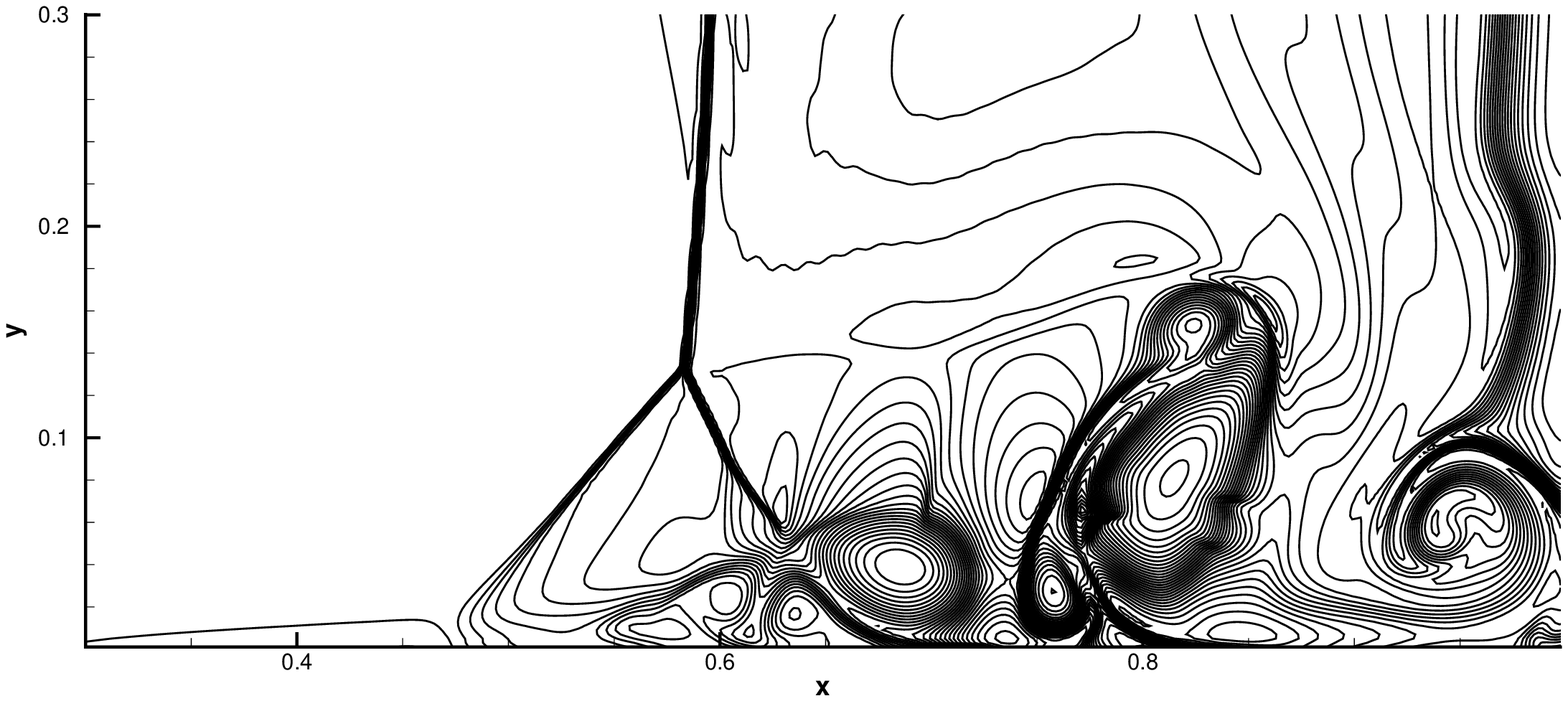}{b}\\
\includegraphics[width=0.55\textwidth]{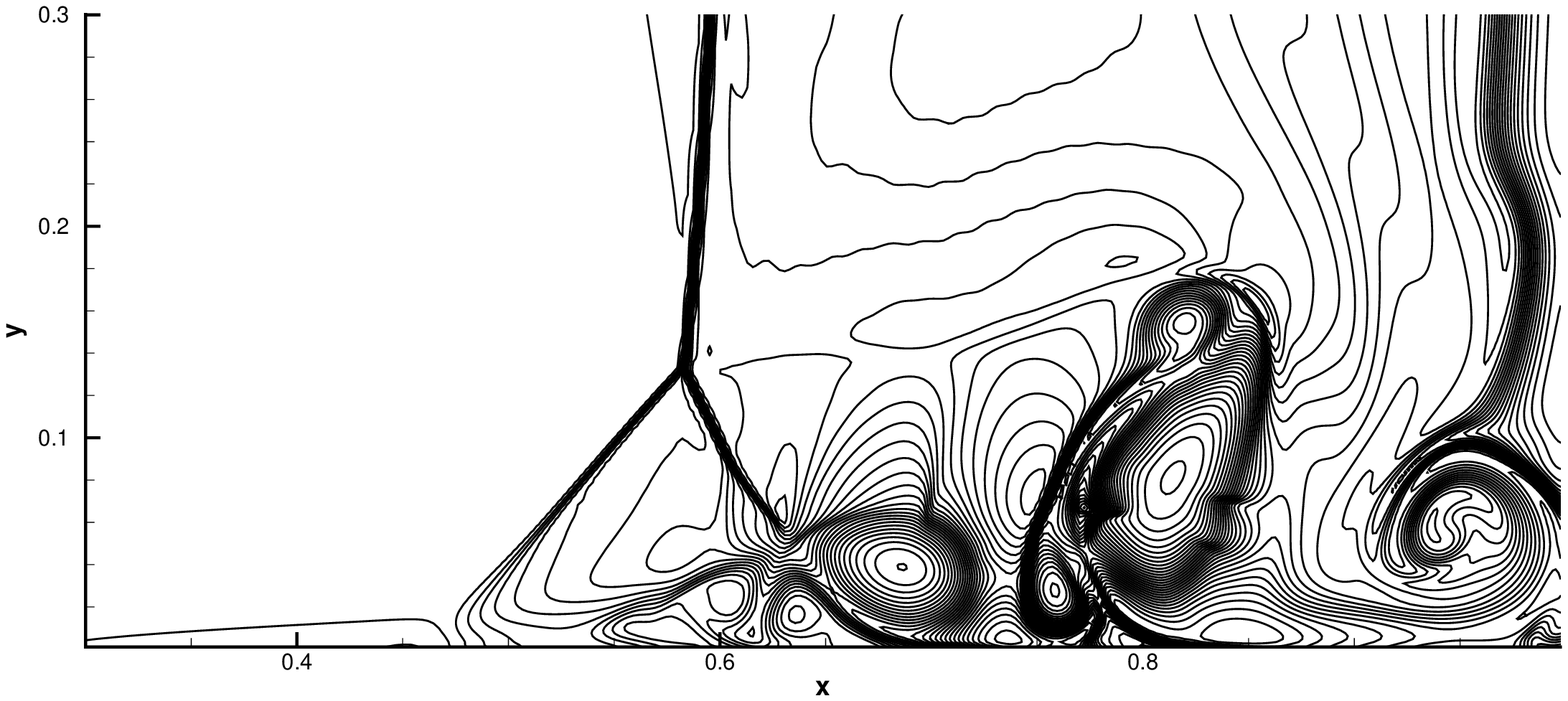}{c}
\caption{\label{vis-shock-2} Viscous shock tube: density
distribution on nonuniform orthogonal mesh (a), nonorthogonal mesh
(b) and uniform mesh (c) with $500\times250$ cells.}
\end{figure}

\begin{figure}[!h]
\centering
\includegraphics[width=0.7\textwidth]{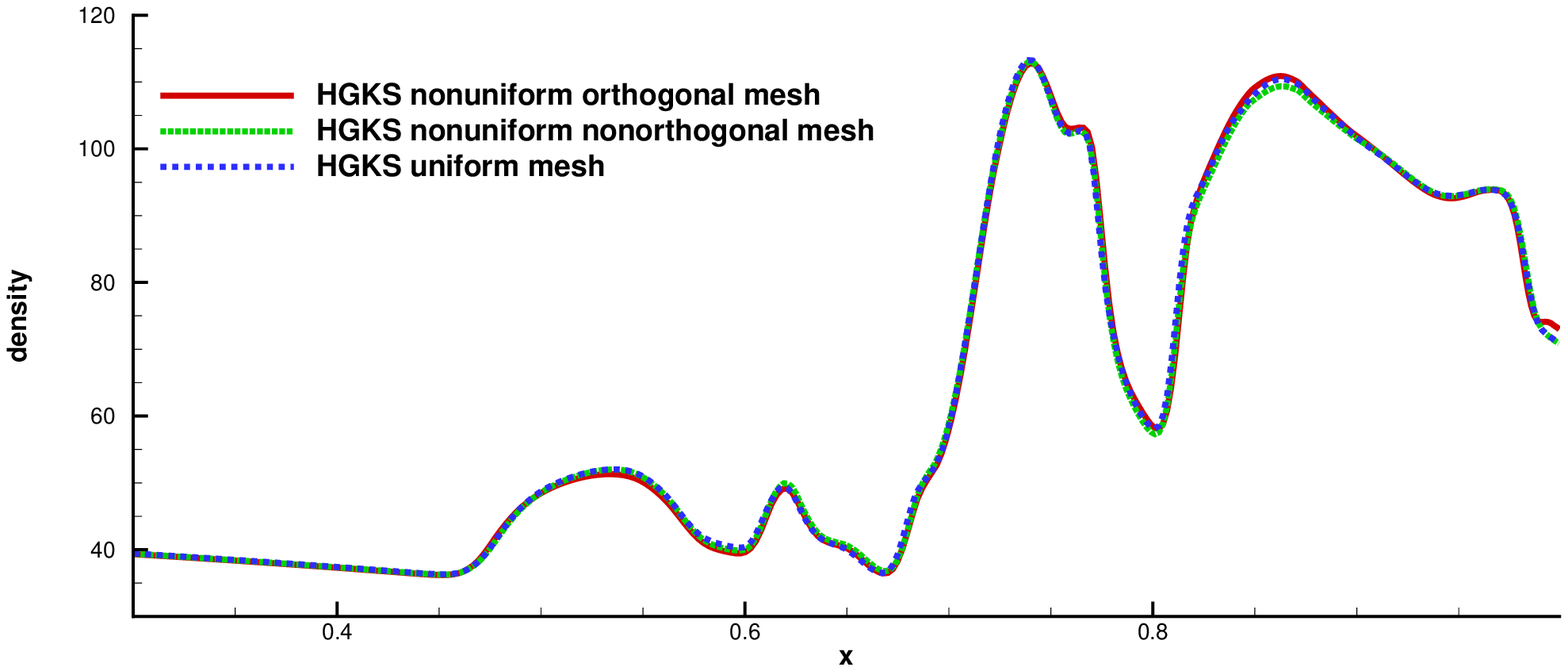}
\caption{\label{vis-shock-3} Viscous shock tube: The density
profiles along the lower wall for different meshes.}
\end{figure}

The membrane is removed at time zero and wave interaction occurs. A
shock wave, followed by a contact discontinuity, moves to the right
with Mach number $\mbox{Ma} =2.37$ and reflects at the right end wall. After
the reflection, it interacts with the contact discontinuity. The
contact discontinuity and shock wave interact with the horizontal
wall and create a thin boundary layer during their propagation. The
solution will develop complex two-dimensional
shock/shear/boundary-layer interactions. The density distributions
on the orthogonal and nonorthogonal meshes with $500\times250$ cells
are presented in Fig.\ref{vis-shock-2}. As reference, the density
distributions on orthogonal uniform mesh with $500\times250$ cells
are presented in Fig.\ref{vis-shock-2} as well. The results match
well with each other.  The density profiles along the lower wall for
$\mbox{Re}=200$ are also presented in Fig.\ref{vis-shock-3}, and numerical
results deviate with other sightly due to different mesh size along
the lower wall.

\section{Conclusion}
In this paper, a two-stage fourth-order gas-kinetic scheme in
curvilinear coordinates is developed for the Euler and Navier-Stokes
solutions. With the two-stage temporal discretization
\cite{GRP-high-1,S2O4-GKS-1}, a reliable framework is provided for
constructing a fourth-order scheme under the gas-kinetic framework.
More importantly, this scheme is as robust as the second-order
scheme and works perfectly for complicated flow simulation. To treat
practical problems with general geometry, such as the turbulent
boundary layer and the flow over a wing-body configuration, the
development of a three-dimensional HGKS in general curvilinear
coordinates becomes necessary. To achieve the high-order accuracy,
the dimension-by-dimension WENO-type reconstruction is adopted in
the computational domain, where the reconstructed Jacobian and the
product of flow variables and local Jacobian are used to get the
point-wise values and spatial derivatives of conservative variables
at Gaussian quadrature points in the computational domain. However,
for the gas-kinetic flow solver, the spatial derivatives of
conservative variables in the physical domain is needed as well,
which is obtained through a procedure of orthogonalization and chain
rule in the local orthogonal coordinates for the flux evaluation in
the normal direction. A variety of numerical tests from the accuracy
test to the solutions with strong discontinuities are presented to
validate the accuracy and robustness of the current scheme. The
geometrical conservation law is precisely satisfied by the current
scheme as well. The current development of HGKS provides a valuable
high-order method for the complicated flow simulation in the complex
geometry under non-uniform non-orthogonal meshes.

\section*{Ackonwledgement}
The current research of L. Pan is supported by National Science Foundation of
China (11701038) and the Fundamental Research Funds for the Central
Universities. The work of K. Xu is supported by
National Science Foundation of China (11772281, 91852114) and Hong Kong research grant council (16206617).

\end{document}